\documentclass[11pt]{amsart}

\usepackage[margin=1in]{geometry}
\usepackage{amsmath,amssymb,amsthm,mathtools, xcolor}
\usepackage{graphicx}
\usepackage{tikz}
\usetikzlibrary{arrows.meta}
\usepackage{placeins}
\usepackage[hidelinks]{hyperref}
\usepackage{mathdots}
\usepackage{booktabs}
\usepackage{longtable}

\newtheorem{theorem}{Theorem}[section]
\newtheorem{proposition}[theorem]{Proposition}
\newtheorem{lemma}[theorem]{Lemma}
\newtheorem{corollary}[theorem]{Corollary}
\newtheorem{remark}[theorem]{Remark}

\newcommand{\Fplus}[2]{F_{#1;#2}}
\newcommand{\Fnotplus}[3]{F_{#1,#2;#3}}
\newcommand{\widehatFplus}[2]{\widehat F_{#1;#2}}
\newcommand{\rtwotilde}{\widetilde p}
\newcommand{\rthreetilde}{\widetilde q}

\title{The lower boundary of the  $\{K_2,K_{p},K_{q}\}$-profile}
\author{Grigoriy Blekherman}
\address{School of Mathematics, Georgia Institute of Technology,
686 Cherry Street
Atlanta, GA 30332}\email{greg@math.gatech.edu}

\author{Annie Raymond}
\address{Department of Mathematics and Statistics,
Lederle Graduate Research Tower,
University of Massachusetts Amherst
710 N. Pleasant Street
Amherst, MA 01003} \email{annieraymond@umass.edu}

\thanks{Grigoriy Blekherman was partially supported by the Simons Travel Award for Mathematicians SFI-MPS-TSM-00023534. Annie Raymond was partially supported by NSF grant DMS-2338532. Some of this material is based upon work supported by the National Science Foundation under Grant No. DMS-1928930, while the second author was in residence at the Simons Laufer Mathematical Sciences Institute in Berkeley, California, during the Spring 2025 semester.}

\begin{document}

\maketitle

\begin{abstract}
We determine the asymptotically sharp minimum number of $q$-cliques that appear in any $n$-vertex graph with $\gamma n^2$ edges and $\delta n^{p}$ $p$-cliques for $q>p>2$. The question of determining the minimum number of $r$-cliques that appear in any $n$-vertex graph with $\gamma n^2$ edges was previously resolved by Razborov for $r=3$, Nikiforov for $r=4$ and Reiher for all $r\geq 5$. For all $r\geq 3$,  the minimum is attained by a complete multipartite graph where all parts are equal except for one that can be smaller or by a limit of such graphs. Very recently, Ma, Wang and Zhu as well as Li, Liu and Zhang showed that the same construction is also a minimizer for the number of $r$-cliques when the number of $s$-cliques is fixed for any $2\leq s<r$. A minimizer in our setting is the disjoint union of the same complete multipartite graph with an independent set. This establishes the lower boundary of the three-dimensional $\{K_2, K_p, K_q\}$-profile.
\end{abstract}

\section{Introduction}

Let $t(H,G)$ be the homomorphism density from a graph $H$ to a target graph $G$, i.e., the probability that a uniformly random (not necessarily injective) map from $V(H)$ to $V(G)$ sends edges of $H$ to edges of $G$. For a finite family $\mathcal U=\{H_1,\ldots,H_\ell\}$, define its
\emph{(rescaled) graph profile} by
$$
\operatorname{cl}\left\{
\left(\frac{t(H_1,G)}{|V(H_1)|!},\ldots,\frac{t(H_\ell,G)}{|V(H_\ell)|!}\right):
G\text{ is a finite graph}
\right\}.
$$ 
The usual graph profile is the closure of the points $(t(H_1,G), \ldots, t(H_\ell,G))$, but the factorial rescaling fits our setting better. 
Since $H_i=K_r$, the rescaled quantities correspond to the number of unlabelled $r$-cliques divided by $|V(G)|^r$, which is exactly what is given by the weighted clique densities introduced below.

Graph profiles capture all valid inequalities in $t(H_1, G)$, $\ldots$, $t(H_\ell, G)$, and as such, they are a central object in extremal combinatorics. It is expected that the problem of understanding all graph profiles is extremely difficult. We know that the problem of certifying polynomial inequalities in homomorphism densities is undecidable \cite{IR, HatamiNorin2011, undecidability}. Still, several two-dimensional graph profiles have been fully determined over the years:
 \begin{itemize}
\item edge and path with two edges (Ahlswede and Katona 1978 \cite{AhlswedeKatona})
\item edge and star (Reiher and Wagner 2017 \cite{ReiherWagner} (corrected by Cairncross and Mubayi 2025 \cite{CairncrossMubayi}))
\item edge and path with an even number of edges (Yang and Zeng 2026 \cite{YangZeng})
\item path with two edges and triangle (Glebov, Grzesik, Hu, Hubai, Kr\'{a}l' and Volec 2016 \cite{GGHHKV})
\item edge and $K_3$ (Razborov 2008 \cite{Razborov})
\item edge and $K_4$ (Nikiforov 2011 \cite{Nikiforov})
\item edge and $K_p$ (Reiher 2016 \cite{Reiher}) 
\item $K_p$ and $K_q$ (Ma, Wang, and Zhu 2026 \cite{MWZ} and Li, Liu and Zhang 2026 \cite{LLZ})
\end{itemize}
Stability results are known in the general case of $K_p$ versus $K_q$. A few other profiles can be understood through general principles: edge and $C_{2r}$, edge and a tree $T$ with a perfect matching, edge and $K_{r,r}$. Profiles have also been studied in the induced setting, e.g., edge against $K_{1,2}$ \cite{LMR}, most pairs of graphs on three vertices \cite{HLNPS, GGHHKV}, as well as the upper profile for $\bar{K}_p$ against $K_q$ \cite{HLNPSthreshold}. Inducibility problems can also be viewed as one-dimensional profiles, and here more cases are understood, e.g., complete bipartite graphs \cite{LPSS}, Tur\'an graphs \cite{LMZ}, and many graphs on at most six vertices \cite{EZL, BHLP, BP, BGLLPS}. As far as we know, no (non-trivial) profile in three dimensions has ever been computed; all preceding results determine relations between pairs of densities. The problem of understanding binomial inequalities is simpler and can be dealt with through tropicalization. In such cases, several arbitrarily large classes of graphs (which result in tropicalizations in arbitrarily many dimensions) are understood, e.g., complete graphs, stars, paths \cite{tropicalization, paths}.

The goal of this paper is to prove the following theorem which allows us to determine the lower boundary of
the graph profile for $\{K_2,K_{p},K_{q}\}$ when the
$K_{q}$-coordinate is regarded as the `up' direction. In other words, we understand the bottom half of this three-dimensional profile. The upper half seems to be more complicated, unlike in the two-dimensional cases.

\begin{theorem}\label{thm:main}
Let $q>p>2$, let $s\geq p-1$ be an integer, and let
$\alpha\in[0,1/s]$ and $\sigma\in[0,1]$.  Every graph on $n$ vertices with
exactly
$$
\sigma^2\frac{s}{2(s+1)}(1-\alpha^2)n^2
$$
edges and at least
$$
\frac{\sigma^{p}}{(s+1)^{p}}
\binom{s+1}{p}
(1+\alpha)^{p-1}(1-(p-1)\alpha)n^{p}
$$
cliques of size $p$ contains at least
$$
\frac{\sigma^{q}}{(s+1)^{q}}
\binom{s+1}{q}
(1+\alpha)^{q-1}(1-(q-1)\alpha)n^{q}
$$
cliques of size $q$.  The bound is attained by the disjoint
union of an independent set of size $(1-\sigma)n$ and a complete
$(s+1)$-partite graph in which $s$ parts have size
$$
\frac{\sigma(1+\alpha)}{s+1}n
$$
and the remaining part has size
$$
\frac{\sigma(1-s\alpha)}{s+1}n
$$
whenever the part sizes are integral. In particular, the bound is asymptotically sharp.
\end{theorem}

We refer to the graph with parameters $s$, $\alpha$, and $\sigma$ defined
in the theorem as the \emph{candidate graph for those parameters}, or
simply as the \emph{candidate graph} when the parameters are clear.  Note
that, in the previous theorem, the edge count is required to be exact whereas the number of cliques of size $p$ is given by a lower bound. Indeed, when we fix the number of $p$-cliques, the number of $q$-cliques in our candidate graphs is strictly decreasing in the number of edges on the positive patches (see \eqref{eq:higher-lambda2}),
so an ``at least'' hypothesis on the number of edges would not imply the
stated conclusion.

In Section~\ref{sec:fill-profile}, we show that every feasible point in the 
$\{K_2,K_{p}\}$-profile is either realized by the candidate graph for some
parameters $s$, $\alpha$, and $\sigma$, or belongs to the closure of these
graphs.  These are results that can be easily recovered from existing results, but as they are essential for the rest of the paper, we include them here. 

On the positive lower boundary of the $\{K_2,K_{p}\}$-profile, i.e., on the scallops, we have
$\sigma=1$, and $s$ and $\alpha$ are the parameters appearing in Reiher's
paper \cite{Reiher}.  The $s$-th scallop is parametrized by complete
$(s+1)$-partite graphs in which $s$ parts are equal and the remaining part
may be smaller.  For fixed $s$ and $\alpha$, the parameter $\sigma$ shrinks
every part of the complete multipartite graph by the same factor and puts the
remaining $(1-\sigma)n$ vertices into an independent set.  The edge and $K_{p}$ densities in this new graph are the densities in the original multipartite graph 
multiplied by $\sigma^2$ and $\sigma^{p}$ respectively. Thus, as $\sigma$ decreases from
$1$ to $0$, the point corresponding to the original multipartite graph moves along a curve from the scallop to the origin; we say that the curve \emph{flows} from the scallop point to the origin.

\begin{figure}[ht]
\centering
\includegraphics[width=0.85\textwidth]{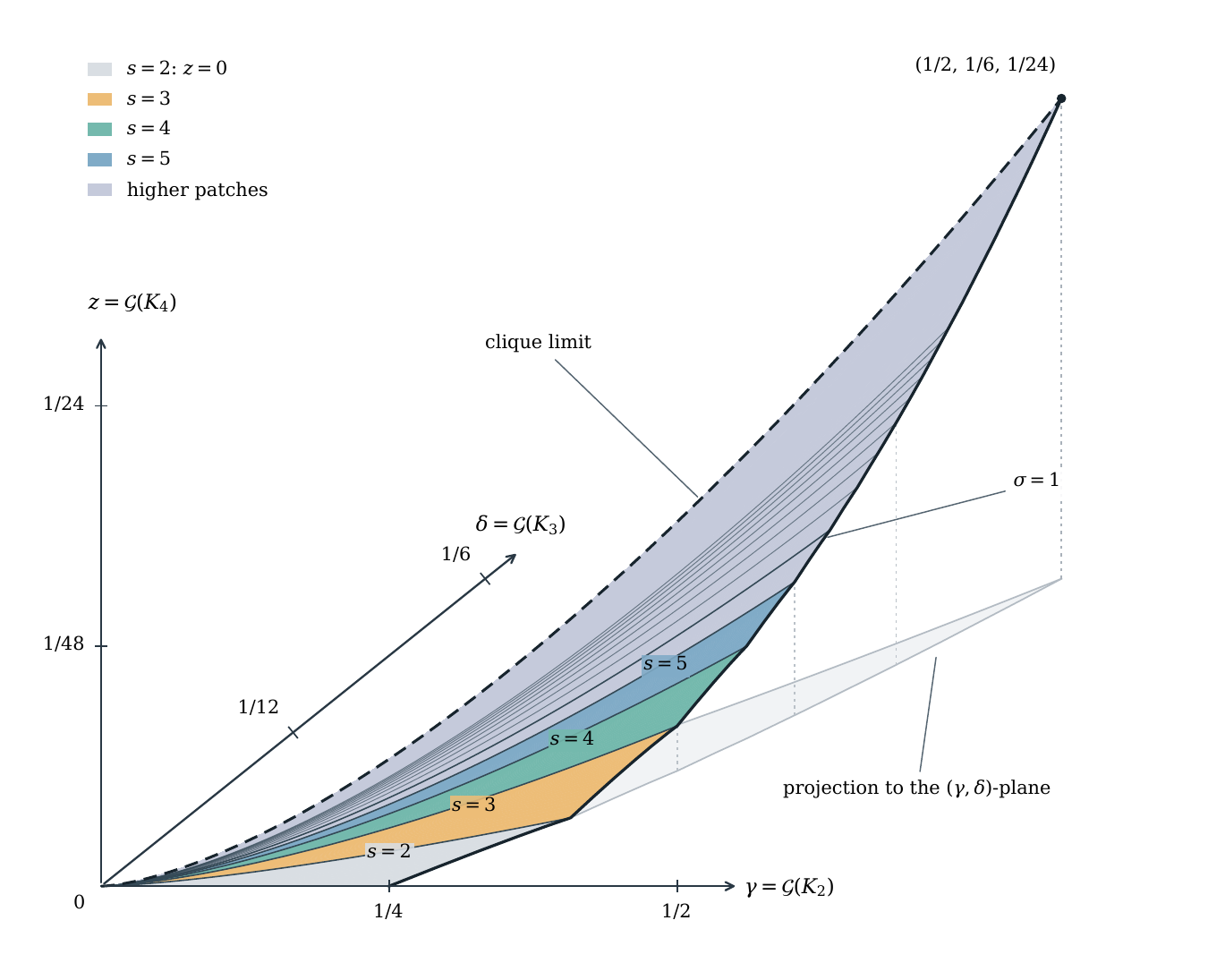}
\caption{The first few patches in the
$\{K_2,K_3,K_4\}$-profile. The projections to the $(\gamma,\delta)$-plane of the heavy solid $\sigma=1$ boundary curves are the scallops (which are concave) and the curves between patches
are the seams (flows from the balanced multipartite graphs to the origin). The dashed curve is the clique-limit curve.}
\label{fig:profile-patches}
\end{figure}

The closure of these flows fills the entire $\{K_2,K_{p}\}$-profile; the upper boundary is obtained as a limit and can be viewed as the flow from the limit point $(\frac{1}{2},\frac{1}{p!})$.  For each $s$, the
scallop and the two flows from its balanced endpoints bound a curvilinear
triangle; for simplicity, we call it the \emph{triangle indexed by $s$}.
Adjacent triangles meet along their common flow from a balanced
multipartite graph, which we call a \emph{seam}.  The triangles lift to
patches on the lower boundary of the
$\{K_2,K_{p},K_{q}\}$-profile.  When $s<q-1$, the corresponding
patches have height zero because their graphs contain no $K_{q}$.  Figure~\ref{fig:profile-patches} illustrates this decomposition when $p=3$ and $q=4$.

We now pass to weighted graphs just like Nikiforov and Reiher did in
\cite{Nikiforov,Reiher}. In particular, throughout the entire paper, we more or less follow the notation from Reiher as well as his strategy, hopefully making this paper slightly easier to read for those already familiar with Reiher's seminal work. Given a
set $X$ and an integer $j\geq0$, let $X^{(j)}$ denote the collection of
all $j$-element subsets of $X$, and put $[n]=\{1,\ldots,n\}$.  A
\emph{weighted graph of order $n$} is a pair
$\mathcal G=(\mathbf{x},\mathbf{a})$, where
$\mathbf{x}=(x_1,\ldots,x_n)$ represents the vertex weights with $x_i\geq0$ and 
$\sum_{i=1}^n x_i=1,$
and
$\mathbf{a}=(a_{ij})_{1\leq i<j\leq n}$ represents the edge weights with $0\leq a_{ij}\leq1.$
For $e=\{i,j\}$, write $a_e=a_{ij}$.  The \emph{$j$-clique density} of
$\mathcal G$ is
$$
\mathcal G(K_j)
=\sum_{M\in[n]^{(j)}}
\prod_{e\in M^{(2)}}a_e\prod_{i\in M}x_i.
$$
In other words, for every set of $j$ vertices in the weighted graph, the weight of that clique is the product of the weights of those $j$ vertices times the product of the weights of the edges between those $j$ vertices. To get $\mathcal{G}(K_j)$, we simply add up the weights of all the $j$-cliques.

The weighted form of the theorem is the statement we shall actually
prove.

\begin{theorem}\label{thm:mainweighted}
Let $q>p>2$, let $s\geq p-1$ be an integer, and let
$\alpha\in[0,1/s]$ and $\sigma\in[0,1]$.  If $\mathcal G$ is a weighted
graph satisfying
$$
\mathcal G(K_2)
=\sigma^2\frac{s}{2(s+1)}(1-\alpha^2)
$$
and
$$
\mathcal G(K_{p})
\geq
\frac{\sigma^{p}}{(s+1)^{p}}
\binom{s+1}{p}
(1+\alpha)^{p-1}(1-(p-1)\alpha),
$$
then
$$
\mathcal G(K_{q})
\geq
\frac{\sigma^{q}}{(s+1)^{q}}
\binom{s+1}{q}
(1+\alpha)^{q-1}(1-(q-1)\alpha).
$$
\end{theorem}

Just as in \cite[Theorem~1.3]{Nikiforov} and \cite[Claim~2.1]{Reiher}, proving Theorem~\ref{thm:mainweighted} yields Theorem~\ref{thm:main}.
\begin{lemma}\label{lem:weighted-implies-finite}
Theorem~\ref{thm:mainweighted} implies Theorem~\ref{thm:main}.
\end{lemma}

\begin{proof}
Let $G=(V,E)$ with $V=[n]$.  Give every vertex weight $1/n$
and put $a_{ij}=1$ when $\{i,j\}\in E$ and $a_{ij}=0$ otherwise.  Then,
for every $j\geq2$, the number of $j$-cliques in $G$ is exactly
$\mathcal G(K_j)n^j$.  The weighted theorem therefore yields the finite
statement.
\end{proof}

Relatedly, passing to weighted graphs does not change the closure of the $\{K_{j_1}, K_{j_2}\}$-profile. Indeed, every finite graph is a weighted graph with uniform vertex weights and edge weights in $\{0,1\}$. Conversely, given a finite weighted graph $\mathcal G=(\mathbf x,\mathbf a)$, form a random graph $G_N$ partitioned into independent sets of sizes $Nx_i+O(1)$, placing each edge between classes $i$ and $k$ independently with probability $a_{ik}$. For each fixed $j\geq2$, the normalized number of $j$-cliques in $G_N$ converges in probability to $\mathcal G(K_j)$: its expectation converges to this value, and its variance tends to zero. Applying this simultaneously to a fixed finite collection of clique orders and choosing suitable realizations gives the required approximation by finite graphs. Thus the ordinary and weighted clique profiles have the same closure.

We next define two functions that we will use throughout the paper. For integers
$2\leq j_1<j_2$ and a feasible value $\gamma$, the first function yields the $K_{j_2}$-density in a complete multipartite graph with equal parts except for one that may be smaller, and which has a $K_{j_1}$ density of $\gamma$. From \cite{MWZ} and \cite{LLZ}, we know this is the minimum among all graphs with a $K_{j_1}$-density of $\gamma$. Define $\Fplus{j_1}{j_2}:[0,\frac{1}{j_1!})\rightarrow [0,\frac{1}{j_2!})$ as follows: given $\gamma\in [0,\frac{1}{j_1!})$, choose the unique integer $s\geq j_1-1$ for which $\gamma\in [\frac{\binom{s}{j_1}}{s^{j_1}},\frac{\binom{s+1}{j_1}}{(s+1)^{j_1}})$ is true, determine the real number $\alpha\in (0,\frac{1}{s}]$ solving $$\gamma=\frac{1}{(s+1)^{j_1}}\binom{s+1}{j_1}(1+\alpha)^{j_1-1}(1-(j_1-1)\alpha)$$ and set

$$\Fplus{j_1}{j_2}(\gamma)=\frac{1}{(s+1)^{j_2}}\binom{s+1}{j_2}(1+\alpha)^{j_2-1}(1-(j_2-1)\alpha).$$

We will often be interested in the case when $j_1=2$ where the formulas specialize to

$$\gamma=\frac{s}{2(s+1)}(1-\alpha^2)=:\gamma_s(\alpha)$$
and

\begin{equation}\label{eq:reiher-profile}
\Fplus{2}{j_2}(\gamma_s(\alpha))=\frac{1}{(s+1)^{j_2}}\binom{s+1}{j_2}(1+\alpha)^{j_2-1}(1-(j_2-1)\alpha).
\end{equation}

As stated before, from \cite{MWZ} and \cite{LLZ}, we know that $\mathcal{G}(K_{j_2})\geq \Fplus{j_1}{j_2}(\mathcal{G}(K_{j_1}))$ over all weighted graphs. Note that $\Fplus{2}{j_2}$ is Reiher's $F_{j_2}$ and $\Fplus{j_1}{j_2}$ is Li--Liu--Zhang's $\Phi_{j_1,j_2}$. We chose to put the size of both cliques involved in the indices because sometimes we will not have $j_1=2$. We are now going to define this function's three-dimensional analog where, again, the number of cliques corresponding to the indices before the semicolon is fixed, and the number of cliques corresponding to the index after the semicolon is minimized.

Moreover, as a convention, we set $\Fplus{2}{2}(x)=x$.

The second function yields the $K_{j_2}$-density in the union of an independent set with a complete multipartite graph with equal parts except for one that may be smaller, and which has an edge density of $\gamma$ and a $K_{j_1}$-density of $\delta$. Let $2<j_1<j_2$, and suppose that $(\gamma,\delta)$ lies in the  $\{K_2,K_{j_1}\}$-profile.  When
$$
\gamma>0
\quad\text{and}\quad
0<\delta<\frac{(2\gamma)^{j_1/2}}{j_1!},
$$
choose $s\geq j_1-1$, $\alpha\in[0,1/s]$, and $\sigma\in(0,1]$ so that
\begin{equation}\label{eq:candidate-parameters}
\gamma=\sigma^2\gamma_s(\alpha),
\qquad
\delta=\sigma^{j_1}\Fplus{2}{j_1}(\gamma_s(\alpha)),
\end{equation}
and define
\begin{equation*}
\Fnotplus{2}{j_1}{j_2}(\gamma,\delta)
=\sigma^{j_2}\Fplus{2}{j_2}(\gamma_s(\alpha)).
\end{equation*}
At the two boundary cases, set
\begin{equation*}
\Fnotplus{2}{j_1}{j_2}(\gamma,0)=0
\end{equation*}
and
\begin{equation*}
\Fnotplus{2}{j_1}{j_2}
\left(\gamma,\frac{(2\gamma)^{j_1/2}}{j_1!}\right)
=\frac{(2\gamma)^{j_2/2}}{j_2!}.
\end{equation*}
The first comes from taking a complete $(j_1-1)$-partite graph with edge density $\gamma$; it has no $j_1$-cliques and thus also no $j_2$-cliques. The second comes from taking a clique with the appropriate edge density (this requires the number of vertices of the graph to go to infinity).

Section~\ref{sec:fill-profile} proves that the parameters exist and that
the definition is unambiguous.  Theorem~\ref{thm:mainweighted} can be rephrased as showing that
\begin{equation*}
\mathcal G(K_{q})
\geq\Fnotplus{2}{p}{q}
\bigl(\mathcal G(K_2),\mathcal G(K_{p})\bigr).
\end{equation*}

Once the theorem is proved, let
$$
\mathcal P=\operatorname{cl}\left\{
\bigl(\mathcal G(K_2),\mathcal G(K_{p}),\mathcal G(K_{q})\bigr):
\mathcal G\text{ is a finite weighted graph}
\right\}.
$$
Then
\begin{equation}\label{eq:minimum}
\Fnotplus{2}{p}{q}(\gamma,\delta)
=\min\{z:(\gamma,\delta,z)\in\mathcal P\}.
\end{equation}

\subsection{Strategy and outline}

\begin{itemize}
\item In Section~\ref{sec:fill-profile}, we first show that the parameters and the
candidate graph for optimality are well defined.  Essentially, the ratio $\delta/\gamma^{p/2}$ chooses the scallop and
the relative size of the smaller part in the complete multipartite component.  The edge density then chooses $\sigma$. Every point of the $\{K_2, K_{p}\}$-profile is associated to a candidate graph or is in the closure of those graphs.

\item In Section~\ref{sec:variational-proof}, we prove
Theorem~\ref{thm:mainweighted} modulo two separation lemmas stated here but proved
later. We broadly justify the use of these separation lemmas in this section too.  The argument for Theorem~\ref{thm:mainweighted} is as follows. Consider
$$
H(\mathcal G):=\mathcal G(K_{q})
-\Fnotplus{2}{p}{q}
\bigl(\mathcal G(K_2),\mathcal G(K_{p})\bigr).
$$
For every edge and $K_{p}$-density in the $\{K_2, K_{p}\}$-profile, the candidate graph associated to that point yields $H(\mathcal{G})=0$. For the sake of contradiction, let us assume that $H(\mathcal{G})$ has a negative minimum. We first choose a counterexample of minimum order and then minimize $H$ over weighted graphs of that order. We would like to use first-order conditions to prove that we have a contradiction, but the candidate function is not differentiable on the seams. So at every point of a positive triangular patch, including on its seams, we construct a smooth local supporting function $$
\widetilde F\leq\Fnotplus{2}{p}{q}
$$
that agrees with $\Fnotplus{2}{p}{q}$ at the point under
consideration.  The counterexample is consequently also a local
minimizer of the smooth function
$$
\widetilde H(\mathcal G)
=
\mathcal G(K_q)
-
\widetilde F\bigl(\mathcal G(K_2),\mathcal G(K_p)\bigr).
$$

First-order (KKT) conditions for $\widetilde H$ give two necessary
inequalities in each \emph{vertex link}, the neighborhood weighted graph as seen from a particular vertex. After normalization, each vertex link gives a point in Ma--Wang--Zhu and  Li--Liu--Zhang's $\{K_{p-1},K_{q-1}\}$-profile, and so we can use their lower bound here. Finally, using the technical lemmas as given, these relationships yield a contradiction. So the goal from then on is to prove the technical lemmas.

\item Section~\ref{sec:preliminaries} goes over some analytic facts about the
$\{K_{p-1},K_{q-1}\}$-profile needed for the separation argument.

\item In Section~\ref{sec:technicallemmas}, we first prove a convex inequality for every possible vertex link with some preliminary work done in Section~\ref{sec:link-comparisons} and then we average that convex inequality and use Jensen's inequality to get a strict inequality, which proves the second technical lemma.
\end{itemize}


\section{Filling the
\texorpdfstring{$\{K_2,K_{p}\}$}{\{K2,Kp\}}-profile}
\label{sec:fill-profile}

The first part of the next monotonicity lemma is \cite[Lemma~3.1]{LLZ}, and the second part follows easily using the same parametrization of the scallop. We reproduce the proof in
our setting for completeness. We first define $$ D_j(N,\beta) := \frac{\binom{N}{j-1}}{j} \frac{N-j+1+j\beta}{(N+\beta)^j}$$ 
for $N\geq j_1-1$ and $0\leq\beta\leq1$, which is nondecreasing as $\beta$ increases.  The pair $$ \bigl(D_{j_1}(N,\beta),D_{j_2}(N,\beta)\bigr) $$ parametrizes the scallop of the $\{K_{j_1},K_{j_2}\}$-profile between the balanced $N$- and $(N+1)$-partite points.

\begin{lemma}[Li--Liu--Zhang, 2026] \label{lem:profile-monotonicities}
Let $2\leq j_1<j_2$, and put
$$x_*=\frac{\binom{j_2-1}{j_1}}{(j_2-1)^{j_1}}.$$
The functions
$$x\longmapsto
\frac{F_{j_1;j_2}(x)}{x^{j_2/j_1}}
\qquad\text{and}\qquad
x\longmapsto
\frac{F_{j_1;j_2}(x)}
{x^{(j_2-1)/(j_1-1)}}
$$
are nondecreasing on $(0,1/j_1!)$ and extend continuously to $x=0$ with value $0$. Both extensions are identically zero on $[0,x_*]$ and strictly increasing on $[x_*,1/j_1!)$.
\end{lemma}

\begin{proof}
On a positive scallop, logarithmic differentiation gives $$ \frac{d}{d\beta} \log\left( \frac{D_{j_2}(N,\beta)} {D_{j_1}(N,\beta)^{j_2/j_1}} \right) = \frac{j_2(j_2-j_1)(1-\beta)} {(N-j_1+1+j_1\beta)(N-j_2+1+j_2\beta)} >0 $$ in the interior, and
\begin{align*}
\frac{d}{d\beta} \log\left( \frac{D_{j_2}(N,\beta)} {D_{j_1}(N,\beta)^{(j_2-1)/(j_1-1)}} \right) ={}& \frac{(j_2-1)(j_2-j_1)(N+1)(1-\beta)} {(N+\beta)(N-j_1+1+j_1\beta)}\\
 &\times \frac{1}{N-j_2+1+j_2\beta} >0. \end{align*}
The function $F_{j_1;j_2}$ vanishes precisely on $[0,x_*]$, whose right endpoint corresponds to the balanced complete $(j_2-1)$-partite graph. Adjacent scallop formulas agree at their common balanced multipartite point. Consequently, the two quotients are continuous and nondecreasing throughout their domain. If $x_*\leq x<y<1/j_1!$, the interval $(x,y)$ contains a nonempty subinterval in the interior of a positive scallop, where both quotients are strictly increasing. Thus both quotients are strictly larger at $y$ than at $x$, proving the asserted global strict monotonicity.
\end{proof}

We also state \cite[Lemma~3.2]{LLZ} as we will need it later; we now allow $\beta$ to be negative in $D_j(N,\beta)$. No corresponding graph exists given that there would be a part with negative size, but the value of $D_j(N,\beta)$ in such cases will allow us to extend a function we care about later on.

\begin{lemma}[Li--Liu--Zhang, 2026] \label{lem:compression}
Let $2\leq j_1<j_2$. If $N\geq j_2-1$ and $-\frac{N-j_2+1}{j_2}\leq \beta \leq 0$, then $D_{j_1}(N, \beta)$ is in the domain of $F_{j_1; j_2}$ and $F_{j_1; j_2}(D_{j_1}(N, \beta))\geq D_{j_2}(N,\beta)$. 
\end{lemma}

Note that the parameter $\sigma$ gets cancelled out in the ratio $\frac{\delta}{\gamma^{p/2}}=\frac{\Fplus{2}{p}(\gamma_s(\alpha))}{\gamma_s(\alpha)^{p/2}}$ which, as a corollary to Lemma~\ref{lem:profile-monotonicities}, is monotone. The ratio thus allows us to find $s$ and $\alpha$.  

\begin{corollary}\label{cor:ratio-branches} For each $s\geq p-1$, the function $$ \alpha\longmapsto \frac{\Fplus{2}{p}(\gamma_s(\alpha))} {\gamma_s(\alpha)^{p/2}} =:g_s(\alpha) $$ is strictly decreasing as $\alpha$ increases from $0$ to $1/s$.  Moreover,
$$
g_s(0)
=
g_{s+1}\left(\frac1{s+1}\right),
$$
and, as $s\geq p-1$ and $\alpha\in[0,1/s]$ vary, the combined range is $0\leq g_s(\alpha) < \frac{2^{p/2}}{p!}$.
The upper endpoint $g_s(\alpha)=\frac{2^{p/2}}{p!}$ is attained only in the limit
corresponding to a clique.
\end{corollary}

\begin{proof} Apply Lemma~\ref{lem:profile-monotonicities} with $(j_1,j_2)=(2,p)$. Since $\gamma_s(\alpha)$ decreases as $\alpha$ increases, the direction of monotonicity is reversed. The endpoint
formulas agree because they represent the same balanced multipartite
graph.  Finally, $$ g_{p-1}\left(\frac1{p-1}\right)=0, $$ while direct evaluation gives $$ \lim_{s\to\infty}g_s(0)=\frac{2^{p/2}}{p!}. $$ \end{proof}

It follows that, whenever
$$
\gamma>0
\quad\text{and}\quad
0<\delta<\frac{(2\gamma)^{p/2}}{p!},
$$
and $(\gamma,\delta)$ is not on a seam, there is a unique $s\geq p-1$ and a unique $\alpha\in(0,1/s)$ such that
\begin{equation}\label{eq:ratio}
\frac{\delta}{\gamma^{p/2}}
=\frac{\Fplus{2}{p}(\gamma_s(\alpha))}
{\gamma_s(\alpha)^{p/2}}.
\end{equation}
In this case, we abuse
language and say that $\delta/\gamma^{p/2}$ \emph{lies in the interior
of scallop $s$}.  If $s\geq q-1$, i.e., if the
$K_{q}$-density of the candidate graph is positive, we call the scallop \emph{positive}.  If $(\gamma, \delta)$ lies on a seam, then there are two choices for $(s,\alpha)$, namely $(s,0)$ and $(s+1,\frac{1}{s+1})$ that yield the same ratio. Throughout the proof of our main theorem, we use the second representation;
equivalently, we always impose
$$
0<\alpha\leq\frac1s.
$$
Thus a seam represented by $(s,1/s)$ corresponds to a balanced
complete $s$-partite component.  The interior of the scallop indexed
by $s$ has positive candidate $K_q$-density precisely when
$s\geq q-1$, while a seam represented by $(s,1/s)$ has positive
candidate $K_q$-density precisely when $s\geq q$.

The factor scaling the multipartite component of the graph is then forced by
\begin{equation}\label{eq:sigma}
\sigma=\sqrt{\frac{\gamma}{\gamma_s(\alpha)}}.
\end{equation}
Equations~\eqref{eq:ratio}--\eqref{eq:sigma} are exactly
\eqref{eq:candidate-parameters}. 

As in previous work, we show that $\Fnotplus{2}{p}{q}$ is continuous. 

\begin{lemma}\label{lem:Fnotplus-continuous}
The function $\Fnotplus{2}{p}{q}$ is continuous on the feasible
$(\gamma,\delta)$-region.
\end{lemma}

\begin{proof}
Let $\gamma>0$ and
$0<\delta<(2\gamma)^{p/2}/p!$.  On the scallop indexed by $s$,
\begin{equation*}
\Fnotplus{2}{p}{q}(\gamma,\delta)
=\gamma^{q/2}
\frac{\Fplus{2}{q}(\gamma_s(\alpha))}
{\gamma_s(\alpha)^{q/2}}.
\end{equation*}
The quotient is continuous within a scallop and adjacent scallop formulas
agree at their common balanced candidate graph.   As $\delta/\gamma^{p/2}$ tends to $2^{p/2}/p!$,  the
scallop index goes to infinity and the quotient goes to 
$2^{q/2}/q!$, which is the value for the endpoint corresponding to a clique.  At the other end, the limiting candidate has a complete
 $(p-1)$-partite component and hence has zero
$K_{q}$-density.

Finally, the standard upper clique-density inequality gives
$$
0\leq
\Fplus{2}{q}(\gamma_s(\alpha))
\leq
\frac{\bigl(2\gamma_s(\alpha)\bigr)^{q/2}}{q!}.
$$
Since $\gamma=\sigma^2\gamma_s(\alpha)$, it follows that
$$
\begin{aligned}
0
\leq
\Fnotplus{2}{p}{q}(\gamma,\delta)
&=
\sigma^q\Fplus{2}{q}(\gamma_s(\alpha))
\\
&\leq
\sigma^q
\frac{\bigl(2\gamma_s(\alpha)\bigr)^{q/2}}{q!}
=
\frac{(2\gamma)^{q/2}}{q!}.
\end{aligned}
$$
This also proves continuity at $(\gamma,\delta)=(0,0)$.
\end{proof}

We now show that Reiher's clique density theorem ensures that the factor $\sigma$, which sets the size of the multipartite component in the candidate graph, is at most one.

\begin{lemma}\label{lem:active-mass}
If $(\gamma,\delta)$ is feasible with $\gamma >0$ and $0<\delta <\frac{(2\gamma)^{p/2}}{p!}$ and $s,\alpha,\sigma$ are chosen by
\eqref{eq:ratio}--\eqref{eq:sigma}, then
\begin{equation*}
\sigma\leq1.
\end{equation*}
\end{lemma}

\begin{proof}
Reiher's clique density theorem gives
$$
\delta\geq\Fplus{2}{p}(\gamma).
$$
On the other hand, \eqref{eq:ratio} gives $$ \frac{\delta}{\gamma^{p/2}} = \frac{\Fplus{2}{p}(\gamma_s(\alpha))} {\gamma_s(\alpha)^{p/2}}. $$ The first monotonicity in Lemma~\ref{lem:profile-monotonicities}, applied with $(j_1,j_2)=(2,p)$, yields that $\frac{\Fplus{2}{p}(x)}{x^{p/2}}$ is strictly increasing throughout its positive range, and therefore gives $\gamma\leq\gamma_s(\alpha)$. By
\eqref{eq:sigma}, this is precisely $\sigma^2\leq1$.
\end{proof}

Corollary \ref{cor:ratio-branches} and Lemma \ref{lem:active-mass} show that every feasible point in the
$\{K_2,K_{p}\}$-profile lies in the closure of the union of the triangles described in the
introduction.  We now prove the monotonicity that allows the
$K_{p}$-hypothesis in Theorems~\ref{thm:main} and
\ref{thm:mainweighted} to be an inequality.

\begin{lemma}\label{lem:Fnotplus-monotone}
For fixed $\gamma$, the function
$$
\delta\longmapsto
\Fnotplus{2}{p}{q}(\gamma,\delta)
$$
is nondecreasing on its feasible interval and strictly increasing on every
scallop on which the candidate $K_{q}$-density is positive.
\end{lemma}

\begin{proof}
On a positive scallop, the $\sigma$-factors cancel in both $\delta/\gamma^{p/2}$ and
$\Fnotplus{2}{p}{q}(\gamma,\delta)/\gamma^{q/2}$, i.e., they yield $\frac{\Fplus{2}{p}(\gamma_s(\alpha))}{\gamma_s(\alpha)^{p/2}}$ and $\frac{\Fplus{2}{q}(\gamma_s(\alpha))}{\gamma_s(\alpha)^{q/2}}$ respectively. 
The logarithmic derivative with respect to $\beta$ computed in
Lemma~\ref{lem:profile-monotonicities} (applied with $p$ and with $q$) is
positive in both cases. After changing to the $\alpha$-parametrization, the two ratios move in the same
direction as $\alpha$ varies.  Adjacent scallops agree at their common
endpoints, and the candidate is identically zero on the scallops indexed by $s< q-1$.
\end{proof}

We now show that Theorems~\ref{thm:main} and
\ref{thm:mainweighted} are sharp since every point of the $\{K_2, K_{p}\}$-profile is either attained by a candidate graph with the right parameters or belongs to the closure of such points. 

\begin{proposition}\label{prop:sharpness}
For every point $(\gamma,\delta)$ in the $\{K_2, K_p\}$-profile, the point 
$(\gamma,\delta,\Fnotplus{2}{p}{q}(\gamma,\delta))$ lies in the closure of the
weighted-graph profile.
If $(\gamma,\delta)$ is not on the upper boundary of the $\{K_2, K_p\}$-profile, it is attained by the finite weighted-graph
analogue of the construction in Theorem~\ref{thm:main}, i.e., by an isolated vertex and a clique where all edges have weight 1 and all vertices have the same weight except for one which might have a smaller weight.
\end{proposition}

\begin{proof}
Suppose first that $\delta=0$.  Reiher's clique density theorem
implies
$$
0\leq\gamma\leq\frac{p-2}{2(p-1)}.
$$
If $\gamma=0$, take an independent weighted graph.  Otherwise, put
$$
\sigma
=
\sqrt{\frac{2(p-1)\gamma}{p-2}}.
$$
Take a balanced complete $(p-1)$-partite weighted graph of total
vertex weight $\sigma$, together with one isolated vertex of weight
$1-\sigma$.  Thus each of the $p-1$ active parts has weight
$\sigma/(p-1)$.  Its edge density is
$$
\binom{p-1}{2}
\left(\frac{\sigma}{p-1}\right)^2
=
\gamma,
$$
and it contains neither a $K_p$ nor a $K_q$.  It therefore realizes
$$
\bigl(\gamma,0,\Fnotplus{2}{p}{q}(\gamma,0)\bigr)
=
(\gamma,0,0).
$$

If $\delta=\frac{(2\gamma)^{p/2}}{p!}$, take a clique of size $m$ where each vertex has weight $\frac{\sqrt{2\gamma}}{m}$ together with an isolated vertex of weight $1-\sqrt{2\gamma}$. Then the $K_j$-density of that graph is $\binom{m}{j}\left(\frac{\sqrt{2\gamma}}{m}\right)^j$ which goes to $\frac{(\sqrt{2\gamma})^j}{j!}$ as $m\rightarrow \infty$.

We may therefore assume that
$\gamma>0$ and $0<\delta<(2\gamma)^{p/2}/p!$.  By
Lemma~\ref{lem:active-mass}, we have $\sigma\leq 1$ so we may form a complete multipartite weighted graph
with $s$ vertices of weight $$\rho:=\frac{\sigma(1+\alpha)}{s+1},$$ one vertex of weight
$\sigma(1-s\alpha)/(s+1)$, and one isolated vertex of weight $1-\sigma$.
Equation~\eqref{eq:reiher-profile} gives
$$
\mathcal G(K_j)=\sigma^j\Fplus{2}{j}(\gamma_s(\alpha)),
\qquad j=2,p,q.
$$
The three coordinates are therefore
$\gamma$, $\delta$, and
$\Fnotplus{2}{p}{q}(\gamma,\delta)$.
\end{proof}


\section{The proof of Theorem~\ref{thm:mainweighted}}
\label{sec:variational-proof}

The proof of the main theorem is given here modulo two convex separation lemmas that are proved in Section~\ref{sec:technicallemmas}. We state the technical lemmas right after we need them for clarity and continuity purposes, and we also explain the general idea behind them.

\subsection{A few facts about weighted graphs}

We first recall some facts about weighted graphs from \cite{Reiher}.  Let
$\mathcal G=(\mathbf{x},\mathbf{a})$ be a finite weighted graph.  For
distinct vertices $i_1,\ldots,i_\ell$, let
$\mathcal G_{i_1\cdots i_\ell}(K_j)$ denote the unlabelled
$K_j$-weight in their common weighted link, i.e., rooted clique weights.  More explicitly, if
$I=[n]\setminus\{i_1,\ldots,i_\ell\}$, then
$$
\mathcal G_{i_1\cdots i_\ell}(K_j)
=
\sum_{M\in I^{(j)}}
\prod_{(b,k)\in[\ell]\times M}a_{i_bk}
\prod_{e\in M^{(2)}}a_e
\prod_{z\in M}x_z.
$$
We shall repeatedly use
\begin{equation}\label{eq:root-average}
\sum_i x_i\mathcal G_i(K_{j-1})=j\mathcal G(K_j)
\end{equation}
and
\begin{equation*}
\frac{\partial\mathcal G(K_j)}{\partial x_i}
=\mathcal G_i(K_{j-1}),
\qquad
\frac{\partial\mathcal G(K_j)}{\partial a_{ik}}
=x_ix_k\mathcal G_{ik}(K_{j-2}).
\end{equation*}
Note that $\mathcal{G}(K_1)=\sum_{i}x_i=1$ and $\mathcal{G}_{ik}(K_0)=1$. 

It is well known that finite weighted graphs cannot lie on the upper boundary of $\{K_2, K_j\}$ except at the origin. We include a proof here for completeness.
\begin{lemma}\label{lem:finiteweightedgraphnotonupperboundary}
 Let $\mathcal G$ be a finite weighted graph with $\gamma=\mathcal G(K_2)>0$. For every integer $j\geq3$,

$$j!\,\mathcal G(K_j)<(2\gamma)^{j/2}.$$
\end{lemma}

\begin{proof}
 Set $A_{ii}=0$ and $A_{ik}=a_{ik}\sqrt{x_ix_k}$ for $i\neq k$, and let $\theta_1,\ldots,\theta_n$ be the eigenvalues of $A$. Expanding $\operatorname{tr}(A^j)$ as a sum over closed walks gives

$$j!\mathcal G(K_j)\leq\operatorname{tr}(A^j),$$
since the product of all edge weights of an ordered clique is at most the product along its corresponding cycle, and the remaining walk terms are nonnegative. The matrix $A$ is nonzero and has trace zero, so at least two of its eigenvalues are nonzero. As $j>2$, it follows that

$$j!\mathcal G(K_j)
\leq\sum_i\theta_i^j
\leq\sum_i|\theta_i|^j
<\left(\sum_i\theta_i^2\right)^{j/2}=\left(2\sum_{i<k}a_{ik}^2x_ix_k\right)^{j/2}
\leq(2\gamma)^{j/2}.$$

\end{proof}

We have the following problem:
$$
\min\left\{
\mathcal G(K_{q})-\Fnotplus{2}{p}{q}(\gamma',\delta') \mid 
\mathcal G(K_1)=1,\ 
\mathcal G(K_2)=\gamma',\ 
\mathcal G(K_{p})= \delta',\
x_i \geq 0 \ \forall i,\
0\leq a_{ij}\leq 1 \ \forall i<j
\right\}.
$$
(The fact that this can be relaxed so that $\mathcal G(K_{p})\geq  \delta'$ can be recovered from Lemma~\ref{lem:Fnotplus-monotone}.)

Using first-order (KKT) optimality conditions, we want to show that 
$\mathcal G(K_{q})
-\Fnotplus{2}{p}{q}
\bigl(\mathcal G(K_2),\mathcal G(K_{p})\bigr)
$
has no negative local minimum at a finite weighted graph
$\mathcal G$ whose vertex weights are positive. However, $\mathcal G(K_{q})
-\Fnotplus{2}{p}{q}
\bigl(\mathcal G(K_2),\mathcal G(K_{p})\bigr)
$ is not differentiable on seams. One possible approach is to minimize $\mathcal G(K_q)$ above because that objective function is differentiable; however, in that case, we would not be able to compute the Lagrange multipliers corresponding to $\mathcal G(K_2)=\gamma$ and $\mathcal G(K_{p})\geq \delta$. A better approach is to introduce a smooth continuation of each triangular patch of $\Fnotplus{2}{p}{q}$ corresponding to letting the smaller part of the multipartite component in the candidate graph have negative size. Outside of the triangular patch at hand, this continuation lies below the adjacent triangular patch, which means that a local minimizer of the problem above will also be a local minimizer for our new smooth function.

For integers $j\geq2$ and $s\geq j-1$, put
$$
f_{j,s}(\alpha)
:=
\frac{1}{(s+1)^j}\binom{s+1}{j}
(1+\alpha)^{j-1}(1-(j-1)\alpha).
$$
For $0\leq\alpha\leq1/s$, we have
$$
f_{j,s}(\alpha)
=
\Fplus{2}{j}(\gamma_s(\alpha)).
$$
Unlike the expression on the right, however,
$f_{j,s}(\alpha)$ continues to describe the fixed $s$-branch when
$\alpha$ is allowed to be slightly larger than $1/s$.

\begin{lemma}\label{lem:continuation}
Let $(\gamma,\delta)$ be a feasible point satisfying
$$
\Fnotplus{2}{p}{q}(\gamma,\delta)>0,
$$
and suppose that it does not lie on the upper clique boundary.  Choose
its representation $(s,\alpha,\sigma)$ so that
$$
0<\alpha\leq\frac1s.
$$
Then there is an open neighborhood $U\subseteq\mathbb R^2$ of
$(\gamma,\delta)$ and a smooth function
$\widetilde F:U\to\mathbb R$ such that
$$
\widetilde F(\gamma',\delta')
\leq
\Fnotplus{2}{p}{q}(\gamma',\delta')
$$
for every feasible $(\gamma',\delta')\in U$,
$$
\widetilde F(\gamma,\delta)
=
\Fnotplus{2}{p}{q}(\gamma,\delta),
$$
and the derivatives of $\widetilde F$ at $(\gamma,\delta)$ are those
obtained from the fixed $s$-branch
$$
(\sigma',\alpha')
\longmapsto
(\sigma')^qf_{q,s}(\alpha').
$$
\end{lemma}

\begin{proof}
Put
$$
R_s(\alpha)
:=
\frac{f_{p,s}(\alpha)}
{\gamma_s(\alpha)^{p/2}}.
$$
At the given point,
$$
\frac{\delta}{\gamma^{p/2}}=R_s(\alpha).
$$
Direct logarithmic differentiation gives
$$
\frac{d}{d\alpha}\log R_s(\alpha)
=
-\frac{p(p-2)\alpha}
{(1+\alpha)(1-\alpha)(1-(p-1)\alpha)}.
$$
Since $\alpha>0$ and the candidate $K_q$-density is positive, this
derivative is nonzero.  The inverse function theorem therefore gives,
on a full neighborhood $U$ of $(\gamma,\delta)$, smooth functions
$\alpha'=\alpha'(\gamma',\delta')$ and
$\sigma'=\sigma'(\gamma',\delta')$ determined by
$$
\frac{\delta'}{(\gamma')^{p/2}}
=
R_s(\alpha'),
\qquad
\sigma'
=
\sqrt{\frac{\gamma'}{\gamma_s(\alpha')}}.
$$
After shrinking $U$ if necessary, define
$$
\widetilde F(\gamma',\delta')
=
(\sigma')^qf_{q,s}(\alpha').
$$
This is smooth on $U$, and its derivatives at $(\gamma,\delta)$ are
precisely those of the fixed $s$-branch.

If
$$
0<\alpha<\frac1s,
$$
we may shrink $U$ so that $0<\alpha'<1/s$ throughout $U$.  At every
feasible point of this neighborhood, the candidate belongs to the
triangle indexed by $s$, and hence
$$
\widetilde F(\gamma',\delta')
=
\Fnotplus{2}{p}{q}(\gamma',\delta').
$$

It remains to consider
$$
\alpha=\frac1s.
$$
This point lies on the seam obtained by flowing from a balanced
complete $s$-partite graph.  Since its candidate $K_q$-density is
positive, we have $s\geq q$.

For $\gamma'$ near $\gamma$, write
$$
\delta_0(\gamma')
:=
(\gamma')^{p/2}R_s\left(\frac1s\right).
$$
Thus
$$
\bigl(\gamma',\delta_0(\gamma')\bigr)
$$
lies on the nearby portion of the same seam.  On the side
$\delta'\geq\delta_0(\gamma')$, the candidate belongs to the triangle
indexed by $s$, and therefore
$$
\widetilde F(\gamma',\delta')
=
\Fnotplus{2}{p}{q}(\gamma',\delta').
$$
On the other side, the candidate belongs to the adjacent triangle
indexed by $s-1$.

Write
$$
\rho=\frac{\sigma}{s},
\qquad
C=\frac{(p-1)!(q-2)}{(p-2)(q-1)!}.
$$
Differentiation at fixed $\gamma$ gives, on the seam,
$$
\frac{\partial\widetilde F}{\partial\delta}
=
C\frac{(s-p+1)!}{(s-q+1)!}\rho^{q-p},
$$
whereas the derivative from the adjacent triangle is
$$
\left.
\frac{\partial\Fnotplus{2}{p}{q}}{\partial\delta}
\right|_{s-1}
=
C\frac{(s-p)!}{(s-q)!}\rho^{q-p}.
$$
Since $s\geq q$, the first quantity is strictly larger than the
second.  This strict inequality persists in a sufficiently small
neighborhood of the seam.

Let $F_{s-1}$ denote the candidate function on the adjacent side and
put
$$
J(\gamma',\delta')
=
F_{s-1}(\gamma',\delta')
-
\widetilde F(\gamma',\delta').
$$
Along the seam,
$$
J\bigl(\gamma',\delta_0(\gamma')\bigr)=0,
$$
whereas throughout the adjacent side of a sufficiently small
neighborhood,
$$
\frac{\partial J}{\partial\delta}<0.
$$
Thus, for fixed $\gamma'$, the function $J(\gamma',\delta')$
decreases to zero as $\delta'$ increases to
$\delta_0(\gamma')$.  Consequently,
$$
J(\gamma',\delta')\geq0
$$
whenever $\delta'\leq\delta_0(\gamma')$.  Hence
$$
\widetilde F
\leq
\Fnotplus{2}{p}{q}
$$
at every feasible point of $U$, with equality at
$(\gamma,\delta)$.
\end{proof}

\begin{proposition}
\label{prop:smooth}
Let $\mathcal G$ be a finite weighted graph whose vertex weights are
positive, and suppose that
$$
\bigl(\mathcal G(K_2),\mathcal G(K_p)\bigr)
$$
does not lie on the upper clique boundary at positive edge density.
Then 
$$
H(\mathcal G)
=
\mathcal G(K_{q})
-\Fnotplus{2}{p}{q}
\bigl(\mathcal G(K_2),\mathcal G(K_{p})\bigr)
$$
cannot have a negative local minimum at $\mathcal G$.
\end{proposition}

\begin{proof}
Suppose $H(\mathcal{G})$ has a negative local minimum at a finite weighted graph $\mathcal{G}$ whose vertex weights are positive.

Write
$$
\gamma=\mathcal G(K_2),
\qquad
\delta=\mathcal G(K_{p}).
$$

Since $\mathcal G(K_q)\geq0$ and $H(\mathcal G)<0$, we have
$$\Fnotplus{2}{p}{q}(\gamma,\delta)>0.$$
Thus $(\gamma,\delta)$ lies in a positive triangular patch or on a
positive seam.  Choose its representation $(s,\alpha,\sigma)$ as in
Lemma~\ref{lem:continuation}, so that
$$
0<\alpha\leq\frac1s.
$$
In particular, at a seam we use the representation belonging to the
patch with the larger index.  Positivity gives $s\geq q-1$.

The common weight of the $s$ equal parts in the multipartite component
of the candidate graph is
$$
\rho=\frac{\sigma(1+\alpha)}{s+1}.
$$

The three coordinates of the candidate graph are
$$
\begin{aligned}
\gamma
&=\sigma^2\frac{s}{2(s+1)}(1-\alpha^2),\\
\delta
&=\sigma^{p}\frac{1}{(s+1)^{p}}\binom{s+1}{p}
(1+\alpha)^{p-1}(1-(p-1)\alpha),\\
\Fnotplus{2}{p}{q}(\gamma,\delta)
&=\sigma^{q}\frac{1}{(s+1)^{q}}\binom{s+1}{q}
(1+\alpha)^{q-1}(1-(q-1)\alpha).
\end{aligned}
$$

Let $\widetilde F$ be given by
Lemma~\ref{lem:continuation}, and define near
$\mathcal G$
$$
\widetilde H(\mathcal G')
=
\mathcal G'(K_{q})
-\widetilde F\bigl(
\mathcal G'(K_2),\mathcal G'(K_{p})
\bigr).
$$
The local comparison
$\widetilde F\leq
\Fnotplus{2}{p}{q}$
gives $\widetilde H(\mathcal G')\geq H(\mathcal G')$ for every weighted graph $\mathcal G'$ sufficiently close to $\mathcal G$, while equality holds at $\mathcal G$.
Therefore $\mathcal G$ is also a negative local minimizer of the smooth
function $\widetilde H$.

The parameters $\alpha$ and $\sigma$ provide two independent
directions along the surface of $\widetilde F$ since $\alpha>0$ (at $\alpha=0$, both $\frac{\partial \gamma}{\partial\alpha}$ and $\frac{\partial \delta}{\partial\alpha}$ vanish).  Consequently,
$$
\frac{\partial}{\partial\alpha}
\widetilde F\bigl(\gamma(\alpha,\sigma),
\delta(\alpha,\sigma)\bigr)
=
\frac{\partial\widetilde F}{\partial\gamma}
\frac{\partial\gamma}{\partial\alpha}
+
\frac{\partial\widetilde F}{\partial\delta}
\frac{\partial\delta}{\partial\alpha}
$$
and
\begin{equation}\label{eq:higher-chainrulesigma}
\frac{\partial}{\partial\sigma}
\widetilde F\bigl(\gamma(\alpha,\sigma),
\delta(\alpha,\sigma)\bigr)
=
\frac{\partial\widetilde F}{\partial\gamma}
\frac{\partial\gamma}{\partial\sigma}
+
\frac{\partial\widetilde F}{\partial\delta}
\frac{\partial\delta}{\partial\sigma}.
\end{equation}
All derivatives are evaluated at the chosen point. Solving for $\frac{\partial \widetilde F(\gamma,\delta)}{\partial\gamma}$ and $\frac{\partial \widetilde F(\gamma,\delta)}{\partial\delta}$ in the previous two equations gives
\begin{equation}\label{eq:higher-lambda2}
\lambda_2:=\frac{\partial \widetilde F(\gamma,\delta)}{\partial\gamma}
=
-\frac{(q-p)q\binom{s+1}{q}}{(p-2)(s+1)s}\rho^{q-2}
\end{equation}
and
\begin{equation}\label{eq:higher-lambdam}
\lambda_{p}:=\frac{\partial \widetilde F(\gamma,\delta)}{\partial\delta}
=
\frac{(q-2)q\binom{s+1}{q}}{(p-2)p\binom{s+1}{p}}\rho^{q-p}.
\end{equation}

Note that $$ \frac{\partial\widetilde F(\gamma(\alpha,\sigma),\delta(\alpha,\sigma))}{\partial \sigma}=\frac{q\Fnotplus{2}{p}{q}(\gamma,\delta)}{\sigma}, \qquad \frac{\partial\delta}{\partial \sigma}=\frac{p\delta}{\sigma}, \qquad \frac{\partial\gamma}{\partial \sigma}=\frac{2\gamma}{\sigma}.$$
Plugging these into~\eqref{eq:higher-chainrulesigma} yields
\begin{equation}\label{eq:general-Euler-homogeneity}
2\lambda_2\gamma+p\lambda_{p}\delta
=q\Fnotplus{2}{p}{q}(\gamma,\delta).
\end{equation}
This can also be seen as an application of Euler's theorem for homogeneous functions.

Moreover, the numbers $\lambda_2$ and $\lambda_{p}$ are respectively
the Lagrange multipliers for the edge and $K_{p}$ equations
in the problem
\begin{equation}\label{eq:optproblem}
\min\left\{
\mathcal G'(K_{q})-\widetilde F(\gamma',\delta') \mid 
\mathcal G'(K_1)=1,\ 
\mathcal G'(K_2)=\gamma',\ 
\mathcal G'(K_{p})= \delta', \
x_i\geq 0 \ \forall i,\
0\leq a_{ij} \leq 1 \ \forall i<j
\right\}
\end{equation}
where $\gamma'$ and $\delta'$ are auxiliary variables such that $(\gamma',\delta')$ is restricted to the neighborhood where $\widetilde F$ is defined. On the feasible set, the objective is
exactly
$$
\mathcal G'(K_q)
-
\widetilde F\bigl(\mathcal G'(K_2),\mathcal G'(K_p)\bigr)
=
\widetilde H(\mathcal G').
$$
Since $\mathcal G$ is a local minimizer of
$\widetilde H$, the point
$$
(\mathcal G,\gamma,\delta)
$$
is a local minimizer of this augmented problem.

Note here that the objective function as well as all constraints are differentiable. Moreover, the auxiliary $\gamma'$ and $\delta'$ coordinates make the three equality gradients independent, and the active edge bounds admit an inward direction. The Mangasarian-Fomovitz constraint qualification (MFCQ) holds, and we can therefore use the theory of first-order conditions. Introduce multipliers $\lambda_1,\lambda_2,\lambda_p$ for the three
equality constraints and the usual nonnegative multipliers
$\eta_i,\mu_{ij},\nu_{ij}$ for the box constraints. The Lagrangian is
\begin{align*}
\mathcal L
=&
\mathcal G'(K_{q})-\widetilde F(\gamma',\delta')
-\lambda_{p}\left(\mathcal G'(K_{p})-\delta'\right)
-\lambda_2\left(\mathcal G'(K_2)-\gamma'\right)
-\lambda_1\left(\mathcal G'(K_1)-1\right)\\
&- \sum_i \eta_i x_i - \sum_{i<j} \mu_{ij} a_{ij} - \sum_{i<j} \nu_{ij} (1-a_{ij}).
\end{align*}

Stationarity with respect to the auxiliary variables gives
$$
0
=
\frac{\partial\mathcal L}{\partial\gamma'}
=
-\frac{\partial\widetilde F}{\partial\gamma}
+\lambda_2
$$
and
$$
0
=
\frac{\partial\mathcal L}{\partial\delta'}
=
-\frac{\partial\widetilde F}{\partial\delta}
+\lambda_p.
$$
Therefore
$$
\lambda_2
=
\frac{\partial\widetilde F}{\partial\gamma},
\qquad
\lambda_p
=
\frac{\partial\widetilde F}{\partial\delta}.
$$
Thus the two quantities computed above by differentiating the
candidate surface are precisely the Lagrange multipliers associated
with the $K_2$- and $K_p$-density constraints.

\subsubsection*{Vertex first-order (KKT) conditions}

Since we are assuming our weighted graph $\mathcal{G}$ is optimal for the previous optimization problem \eqref{eq:optproblem}, stationarity with respect to $x_i$ gives us that $\frac{\partial \mathcal L}{\partial x_i}=0$ for every $i$. Moreover, since we are assuming that every vertex weight is positive, complementary slackness gives us that $\eta_i=0$ for every $i$. Thus we have
\begin{equation}\label{eq:general-vertex-equation}
\mathcal G_i(K_{q-1})
-\lambda_{p}\mathcal G_i(K_{p-1})
-\lambda_2\mathcal G_i(K_1)
=\lambda_1
\end{equation}
for every $i$.  Therefore, $$\sum_{i} x_i(\mathcal G_i(K_{q-1})-\lambda_{p}\mathcal G_i(K_{p-1})
-\lambda_2\mathcal G_i(K_1))=\sum_{i} x_i\lambda_1=\lambda_1$$ which, using \eqref{eq:root-average}, is equivalent to
$$
\lambda_1
=
q\mathcal G(K_{q})
-p\lambda_{p}\delta
-2\lambda_2\gamma.
$$
Thus \eqref{eq:general-Euler-homogeneity} shows explicitly that
\begin{equation}\label{eq:general-nu}
\lambda_1
=
q
\bigl(\mathcal G(K_{q})-\Fnotplus{2}{p}{q}(\gamma,\delta)\bigr)
\end{equation}
and so $\lambda_1$ is negative since we assumed that $\mathcal{G}(K_{q})-\Fnotplus{2}{p}{q}(\gamma,\delta)$ is negative.

Substituting~\eqref{eq:general-nu} into~\eqref{eq:general-vertex-equation} and plugging in $\lambda_2$ and $\lambda_{p}$  yields the additional relationship that

\begin{equation}\label{eq:finalformvertexkkt}
\mathcal{G}_i(K_{q-1})=-\frac{(q-p)q\binom{s+1}{q}}{(p-2)(s+1)s}\rho^{q-2}\mathcal{G}_i(K_1)+\frac{(q-2)q\binom{s+1}{q}}{(p-2)p\binom{s+1}{p}}\rho^{q-p}\mathcal{G}_i(K_{p-1})+q(\mathcal{G}(K_{q})-\Fnotplus{2}{p}{q}(\gamma,\delta)).
\end{equation}

\subsubsection*{Edge first-order (KKT) conditions}

Since we are assuming that our weighted graph $\mathcal{G}$ is optimal for our optimization problem, stationarity with respect to $a_{ij}$ yields that $\frac{\partial \mathcal L}{\partial a_{ij}}=0$ for every $i<j$.

We claim that 
\begin{equation*}
a_{ij}(\mathcal G_{ij}(K_{q-2})-\lambda_{p}\mathcal G_{ij}(K_{p-2})-\lambda_2)\leq 0.
\end{equation*}
 If $a_{ij}=0$, this is trivially true.  Otherwise, $a_{ij}>0$, and so by complementary slackness, we know that $\mu_{ij}=0$. Thus $\frac{\partial \mathcal L}{\partial a_{ij}}=0$ yields
 
 $$x_i x_j (\mathcal G_{ij}(K_{q-2})-\lambda_{p}\mathcal G_{ij}(K_{p-2})-\lambda_2) = -\nu_{ij}$$
where $\nu_{ij}\geq 0$, and thus the claim holds.

We are now looking for a different relationship between $\mathcal{G}_i(K_{q-1})$, $\mathcal{G}_i(K_{p-1})$ and $\mathcal G_i(K_1)$ to match the one we had obtained from the vertex first order conditions. Our claim implies that $$\sum_{j\neq i} x_j a_{ij}(\mathcal G_{ij}(K_{q-2})-\lambda_{p}\mathcal G_{ij}(K_{p-2})-\lambda_2) \leq 0$$ which is equivalent to
\begin{equation*}
(q-1)\mathcal G_i(K_{q-1})
-(p-1)\lambda_{p}
\mathcal G_i(K_{p-1})
-\lambda_2\mathcal G_i(K_1)
\leq0.
\end{equation*}
Plugging in $\lambda_2, \lambda_{p}$ in the last inequality as well as $\mathcal{G}_i(K_{q-1})$ as given by \eqref{eq:finalformvertexkkt}, we get 

\begin{align*}
-\frac{(q-2)(q-p)q\binom{s+1}{q}}{(p-2)(s+1)s}\rho^{q-2}\mathcal{G}_i(K_1)+\frac{(q-p)(q-2)q\binom{s+1}{q}}{(p-2)p\binom{s+1}{p}}\rho^{q-p}\mathcal{G}_i(K_{p-1})\\
+q(q-1)(\mathcal{G}(K_{q})-\Fnotplus{2}{p}{q}(\gamma,\delta))\leq 0
\end{align*}
which can be rewritten as
\begin{equation}\label{eq:finaledgekkt}
\frac{(p-1)!\mathcal{G}_i(K_{p-1})}{\rho^{p-1}} \leq \frac{(s-1)!}{(s+1-p)!\rho} \mathcal{G}_i(K_1) + \frac{q!(s+1-q)!(p-2)(q-1)}{(q-2)(q-p)(s+1-p)!\rho^{q-1}}(\Fnotplus{2}{p}{q}(\gamma,\delta)-\mathcal{G}(K_{q}))
\end{equation}

\subsubsection*{Constructing a vertex link graph}

By looking at a smaller weighted graph on $n-1$ vertices, we can come up with yet another relationship between $\mathcal{G}_i(K_{q-1})$, $\mathcal{G}_i(K_{p-1})$ and $\mathcal G_i(K_1)$.

First note that $\mathcal G_i(K_1)>0$ for every $i$.  Indeed, if it
were zero,  then we would also have $\mathcal{G}_i(K_{p-1})=0$  and $\mathcal{G}_i(K_{q-1})=0$ which would contradict \eqref{eq:general-vertex-equation} and
the fact that $\lambda_1<0$ by \eqref{eq:general-nu}.  Thus, we can construct a new weighted graph $\mathcal{G}^{*_i}$ of order $n-1$ with vertex weights $x_j^*=\frac{a_{ij}x_j}{\mathcal{G}_i(K_1)}$ for $j\in [n]\backslash\{i\}$ and with edge weights obtained by restricting $\mathbf{a}$ to $([n]\backslash \{i\})^{(2)}$. This graph is called the \emph{link (graph) of vertex $i$}. Its $K_{p-1}$ and $K_{q-1}$-densities are

$$\mathcal{G}^{*_i}(K_{p-1})=\frac{\mathcal{G}_i(K_{p-1})}{\mathcal{G}_i(K_1)^{p-1}}, \textup{ and}$$

$$\mathcal{G}^{*_i}(K_{q-1})=\frac{\mathcal{G}_i(K_{q-1})}{\mathcal{G}_i(K_1)^{q-1}}.$$

By Ma--Wang--Zhu and Li--Liu--Zhang's higher-order clique density lemma, we know that $\mathcal{G}^{*_i}(K_{q-1})\geq \Fplus{p-1}{q-1}(\mathcal{G}^{*_i}(K_{p-1}))$ which, after plugging in the densities of the link, is equivalent to

$$\mathcal{G}_i(K_{q-1}) \geq \mathcal{G}_i(K_1)^{q-1} \Fplus{p-1}{q-1}\left(\frac{\mathcal{G}_i(K_{p-1})}{\mathcal{G}_i(K_1)^{p-1}}\right).$$

By using \eqref{eq:finalformvertexkkt}, this in turn is equivalent to 

\begin{equation}\label{eq:vertexkktinlink}\begin{aligned}
\frac{\mathcal{G}_i(K_1)^{q-1}}{\rho^{q-1}}&\Fplus{p-1}{q-1}\left(\frac{\mathcal{G}_i(K_{p-1})}{\mathcal{G}_i(K_1)^{p-1}}\right)\\
& \leq  \frac{(q-2)q\binom{s+1}{q}}{(p-2)p\binom{s+1}{p}}\frac{\mathcal{G}_i(K_{p-1})}{\rho^{p-1}}-\frac{(q-p)q\binom{s+1}{q}}{(p-2)(s+1)s}\frac{\mathcal{G}_i(K_1)}{\rho}-\frac{q(\Fnotplus{2}{p}{q}(\gamma,\delta)-\mathcal{G}(K_{q}))}{\rho^{q-1}}.
\end{aligned}\end{equation}

\subsubsection*{Averaging}

Finally, thinking about $\mathbf{x}$ as a probability distribution where we choose $i$ with probability $x_i$, we have that

\begin{equation}\label{eq:averagingD}
\mathbb{E}\left(\frac{\mathcal{G}_i(K_1)}{\rho}\right)=\sum_{i} x_i\frac{\mathcal{G}_i(K_1)}{\rho}=\frac{2\mathcal{G}(K_2)}{\rho}=\frac{2\gamma}{\rho}=\sigma s (1-\alpha),
\end{equation}
and
\begin{equation}\label{eq:averagingA}
\mathbb{E}\left(\frac{(p-1)!\mathcal{G}_i(K_{p-1})}{\rho^{p-1}}\right)=\sum_{i} x_i\left(\frac{(p-1)!\mathcal{G}_i(K_{p-1})}{\rho^{p-1}}\right)=\frac{p!\mathcal{G}(K_{p})}{\rho^{p-1}}=\frac{p!\delta}{\rho^{p-1}}=\sigma \frac{p!}{s+1}\binom{s+1}{p} (1-(p-1)\alpha).
\end{equation}

We now introduce some new quantities and rewrite some of these equations in terms of these quantities to see that we have reached the assumptions of Lemma~\ref{lem:higher-smooth-separation}, and that we thus have a contradiction. 

To make the dependence on the link exponents more transparent, we now introduce the notation
$$
\rtwotilde=p-1,
\qquad
\rthreetilde=q-1.
$$
For a real number $x$ and an integer $j\geq0$, we write
$$
(x)_j=x(x-1)\cdots(x-j+1),
\qquad
(x)_0=1.
$$
It will be more convenient to use
the homomorphism density version of $\Fplus{\rtwotilde}{\rthreetilde}$ extended continuously to its clique endpoint, i.e., we want to find the minimum $t(K_{\rthreetilde}, G)$ such that $t(K_{\rtwotilde}, G)\geq x$ whereas we had previously been minimizing $\frac{t(K_{\rthreetilde},G)}{{\rthreetilde}!}$ such that $\frac{t(K_{\rtwotilde},G)}{{\rtwotilde}!}\geq \frac{x}{\rtwotilde!}$. Define $\widehatFplus{\rtwotilde}{\rthreetilde}$ on $[0,1]$
by
\begin{equation}\label{eq:labelled-Fplus}
\widehatFplus{\rtwotilde}{\rthreetilde}(x)
:=
\begin{cases}
\rthreetilde!\,
\Fplus{\rtwotilde}{\rthreetilde}
\left(\dfrac{x}{\rtwotilde!}\right),
&0\leq x<1,\\[2mm]
1,&x=1.
\end{cases}
\end{equation}

The value at $x=1$ is the continuous limiting value corresponding to
a clique.

Define
\begin{equation*}
D_i=\frac{\mathcal G_i(K_1)}{\rho},
\qquad
P_i=\frac{\rtwotilde!\mathcal G_i(K_{\rtwotilde})}
{\rho^{\rtwotilde}},
\qquad
Q_i=\frac{\rthreetilde!\mathcal G_i(K_{\rthreetilde})}
{\rho^{\rthreetilde}}.
\end{equation*}
The actual normalized link point is
$$
\left(
\frac{P_i}{D_i^{\rtwotilde}},
\frac{Q_i}{D_i^{\rthreetilde}}
\right)
=
\left(
\rtwotilde!\mathcal G^{*_i}(K_{\rtwotilde}),
\rthreetilde!\mathcal G^{*_i}(K_{\rthreetilde})
\right).
$$
Since its first coordinate is the labelled
$K_{\rtwotilde}$-density of a weighted graph, we have
$$
0\leq\frac{P_i}{D_i^{\rtwotilde}}\leq1,
$$
or equivalently,
$$
0\leq P_i\leq D_i^{\rtwotilde}.
$$
Though it would feel more natural to work directly with the normalized link point, we shall nevertheless work with $D_i$ and $P_i$, because their means
are controlled linearly by the rooted identities, whereas the coordinates
of the normalized link points are not.

Put
$$
Y=(\rthreetilde-\rtwotilde)!\binom{s-\rtwotilde}{\rthreetilde-\rtwotilde}=(s-\rtwotilde)_{\rthreetilde-\rtwotilde},
\qquad
K=(\rtwotilde-1)!\binom{s-1}{\rtwotilde-1}=(s-1)_{\rtwotilde-1},
$$
and
\begin{equation*}
c=
\frac{(\rtwotilde-1)(\rthreetilde+1)!
\bigl(\Fnotplus{2}{\rtwotilde+1}{\rthreetilde+1}(\gamma,\delta)-\mathcal G(K_{\rthreetilde+1})\bigr)}
{Y\rho^{\rthreetilde}}>0.
\end{equation*}

With this notation in place, equations~\eqref{eq:vertexkktinlink} and \eqref{eq:finaledgekkt} become

\begin{equation}\label{eq:higher-local-profile}
D_i^{\rthreetilde}
\widehatFplus{\rtwotilde}{\rthreetilde}
\left(\frac{P_i}{D_i^{\rtwotilde}}\right)
\leq
\frac{Y}{\rtwotilde-1}
\bigl((\rthreetilde-1)P_i
-(\rthreetilde-\rtwotilde)KD_i-c\bigr),
\end{equation}

\begin{equation}\label{eq:higher-local-edge}
P_i\leq KD_i+
\frac{\rthreetilde}
{(\rthreetilde-1)(\rthreetilde-\rtwotilde)}c,
\end{equation}

and equations~\eqref{eq:averagingD} and \eqref{eq:averagingA} become
\begin{equation}\label{eq:higher-averaging-D}
\mathbb ED
=
\frac{2\mathcal G(K_2)}{\rho}
=\sigma s(1-\alpha)
\end{equation}
and
\begin{equation}\label{eq:higher-averaging-B}
\mathbb EP
=
\frac{(\rtwotilde+1)!\mathcal G(K_{\rtwotilde+1})}
{\rho^{\rtwotilde}}
=
\sigma \frac{(\rtwotilde+1)!}{s+1}\binom{s+1}{\rtwotilde+1} (1-\rtwotilde\alpha)=\sigma(s)_{\rtwotilde}(1-\rtwotilde\alpha).
\end{equation}
Thus every pair $(D_i,P_i)$ satisfies
\eqref{eq:higher-local-profile}--\eqref{eq:higher-local-edge}, while their
barycenter is prescribed by
\eqref{eq:higher-averaging-D}--\eqref{eq:higher-averaging-B}.  This is
forbidden by Lemma~\ref{lem:higher-smooth-separation} below.  We have
reached a contradiction.
\end{proof}

\subsection{Statements of the technical lemmas and the general ideas behind them}\label{sec:generalideatechnicalstatements}

As in \cite{Reiher,LLZ}, we convert inequalities in vertex links into a
global contradiction by averaging. We first explain the convex separation that Lemma~\ref{lem:higher-pointwise} is designed
to produce. As mentioned above, we work with $(D_i,P_i)$ instead of the actual normalized link points because their means are controlled linearly
by the rooted identities.

Let $\mathcal R_c$ be the set of pairs $(D,P)$ satisfying $D>0$, $D^{\rtwotilde} \geq P \geq 0$ and the
two KKT inequalities
\eqref{eq:higher-local-profile}--\eqref{eq:higher-local-edge}.  Every
$(D_i,P_i)$ belongs to $\mathcal R_c$, so their barycenter must belong to
its convex hull.  We shall show that the prescribed barycenter does not
belong to this convex hull when $c>0$, i.e., if the candidate graph is not optimal and we instead have some negative minimizer.

To understand an appropriate separation, let us first consider $c=0$. In this case we cannot find a contradiction because of our candidate graph construction.  Put
$$
K=(s-1)_{\rtwotilde-1}.
$$
For the isolated vertex of the candidate graph, $(D_i,P_i)=(0,0)$.  For the exceptional vertex,
$$
(D_i,P_i)=(s,sK).
$$
For each of the $s$ regular vertices with weight $\rho$, we have
$$
(D_i,P_i)
=
\left(
\frac{s-\alpha}{1+\alpha},
\frac{K\bigl(s-(s(\rtwotilde-1)+\rtwotilde)\alpha\bigr)}
{1+\alpha}
\right).
$$
The barycenter is
$$
(\overline D,\overline P)
=
(1-\sigma)(0,0)
+(\sigma-s\rho)(s,sK)
+s\rho\left(
\frac{s-\alpha}{1+\alpha},
\frac{K\bigl(s-(s(\rtwotilde-1)+\rtwotilde)\alpha\bigr)}
{1+\alpha}
\right).
$$

Note that in the proof of Proposition~\ref{prop:smooth} where we assumed that $c>0$, we showed that $\mathcal{G}_i(K_1)>0$ for every $i$, so in that case we cannot get $(0,0)$ as a point. Thus, one idea for the appropriate separation is to use the line going through both nonzero candidate link points, namely
$$
\frac{\rtwotilde KD-P}{\rtwotilde-1}=sK.
$$
We call the left-hand side $L(D,P):=\frac{\rtwotilde KD-P}{\rtwotilde-1}$. The prescribed barycenter satisfies
$$
L(\overline D,\overline P)=\sigma sK\leq sK.
$$
We would like to prove $L(D,P)>sK$ pointwise when
$c>0$, but that is too strong.

Instead, put
$$
\widehat c
=\frac{\rtwotilde c}
{(\rtwotilde-1)(\rthreetilde-\rtwotilde)}.
$$
Lemma~\ref{lem:higher-pointwise} proves
$$
L(D,P)\geq
sK\left(
1+\frac{(\rtwotilde-1)\widehat c}{\rtwotilde KD}
\right)^{\rtwotilde/(\rtwotilde-1)}
-\widehat c
=:\psi_c(D).
$$

In the $(D,L)$-plane, every allowed point therefore lies above the graph
of the strictly convex function $\psi_c$, and hence above its tangent at
$D=\overline D$: $$L(D,P)\geq \psi_c(\overline D)+\psi'_c(\overline{D})(D-\overline{D}).$$  This tangent is an affine separator in the original
$(D,P)$-coordinates.  In Lemma~\ref{lem:higher-smooth-separation}, averaging, or equivalently, using Jensen's inequality, makes its linear term vanish.  Using the facts that
$\overline D\leq s$ and $c>0$, an application of the Bernoulli inequality $(1+x)^\eta \geq 1+\eta x$ for any real number $\eta \geq 1$ and any real number $x\geq -1$ yields
$$
L(\overline D,\overline P)\geq\psi_c(\overline D)>sK,
$$
contradicting $L(\overline D,\overline P)=\sigma sK\leq sK$.

Having explained the ideas behind these two technical lemmas, we now state them and delay their proofs---which are technical---so that we can finish proving our main theorem.

\begin{lemma}
\label{lem:higher-pointwise}
Let $s\geq\rthreetilde$, let 
$$Y=(s-\rtwotilde)_{\rthreetilde-\rtwotilde},
\qquad
K=(s-1)_{\rtwotilde-1},$$
 let $c>0$, and suppose that
$D>0$, $0\leq P\leq D^{\rtwotilde}$ and that $(D,P)$ satisfies
the KKT inequalities \eqref{eq:higher-local-profile}--\eqref{eq:higher-local-edge}.
Put
\begin{equation*}
\widehat c
=\frac{\rtwotilde c}
{(\rtwotilde-1)(\rthreetilde-\rtwotilde)}.
\end{equation*}
Then
\begin{equation}\label{eq:higher-pointwise}
\frac{\rtwotilde KD-P}{\rtwotilde-1}
\geq
sK\left(
1+\frac{(\rtwotilde-1)\widehat c}{\rtwotilde KD}
\right)^{\rtwotilde/(\rtwotilde-1)}
-\widehat c.
\end{equation}
\end{lemma}

\begin{lemma}
\label{lem:higher-smooth-separation}
Fix $s\geq\rthreetilde$, $0\leq\alpha\leq1/s$,
$0<\sigma\leq1$, and $c>0$.  Put
$$
Y=(s-\rtwotilde)_{\rthreetilde-\rtwotilde},
\qquad
K=(s-1)_{\rtwotilde-1}.
$$
There is no finitely supported probability distribution on pairs $(D,P)$ satisfying
$$
D>0,
\qquad
0\leq P\leq D^{\rtwotilde},
$$
\begin{equation}\label{eq:higher-separation-local-profile}
D^{\rthreetilde}
\widehatFplus{\rtwotilde}{\rthreetilde}
\left(\frac{P}{D^{\rtwotilde}}\right)
\leq
\frac{Y}{\rtwotilde-1}
\bigl((\rthreetilde-1)P
-(\rthreetilde-\rtwotilde)KD-c\bigr),
\end{equation}
and
\begin{equation}\label{eq:higher-separation-local-edge}
P\leq KD+
\frac{\rthreetilde}
{(\rthreetilde-1)(\rthreetilde-\rtwotilde)}c,
\end{equation}
whose means are
\begin{equation}\label{eq:higher-separation-means}
\overline D=\sigma s(1-\alpha),
\qquad
\overline P=\sigma(s)_{\rtwotilde}(1-\rtwotilde\alpha).
\end{equation}
\end{lemma}

\subsection{Completion of the proof}

\begin{proof}[Proof of Theorem~\ref{thm:mainweighted}]
If we can prove the lower bound, then by Proposition~\ref{prop:sharpness} we know that it is best possible since it can be realized with the candidate graphs. It thus remains to establish the lower
bound.

Suppose that $\mathcal{G}(K_{q})<\Fnotplus{2}{p}{q}(\mathcal{G}(K_2),\mathcal{G}(K_{p}))$ for some finite weighted
graph.  Choose its order $n$ to be the smallest possible among all counterexamples,
and then minimize
$$
H(\mathcal G)
=
\mathcal G(K_q)
-\Fnotplus{2}{p}{q}
\bigl(\mathcal G(K_2),\mathcal G(K_p)\bigr)
$$
among weighted graphs of order $n$.  The minimum exists by compactness
of the vertex-weight simplex and edge-weight cube, together with
Lemma~\ref{lem:Fnotplus-continuous}.  Every vertex weight is positive,
since a zero-weight vertex could be deleted, contradicting the minimality
of $n$.

If it has positive edge density, the minimizer cannot lie in a triangle with
$s\leq q-2$, because the candidate
$K_q$-density is zero there.  It also cannot lie on the
upper clique boundary by Lemma~\ref{lem:finiteweightedgraphnotonupperboundary}. 

We are therefore left with a positive triangle (including its seams); Proposition~\ref{prop:smooth} gives
a contradiction in that case. 
 Hence no negative
minimizer exists, and continuity extends the nonnegativity of $H$ to
the closure of the weighted-graph profile.

This proves $\mathcal{G}(K_{q})\geq \Fnotplus{2}{p}{q}(\mathcal{G}(K_2),\mathcal{G}(K_{p}))$ for an exact
$K_p$-density.  Lemma~\ref{lem:Fnotplus-monotone} allows
that density to be larger than the value in
\eqref{eq:candidate-parameters}, giving the stated form of
Theorem~\ref{thm:mainweighted}.  Lemma~\ref{lem:weighted-implies-finite}
then gives Theorem~\ref{thm:main}.

Proposition~\ref{prop:sharpness} supplies the required
points in the closure.
\end{proof}

\section{Analytic facts about the
\texorpdfstring{$\{K_{\rtwotilde},K_{\rthreetilde}\}$}
{(Kp-1,Kq-1)}-profile}
\label{sec:preliminaries}

Our aim now is to prove
Lemma~\ref{lem:higher-smooth-separation}.  First we establish the analytic
facts about the labelled profile
$\widehatFplus{\rtwotilde}{\rthreetilde}$ introduced in
\eqref{eq:labelled-Fplus}.

We state the monotonicities of Lemma~\ref{lem:profile-monotonicities} in terms of $\widehat{F}$. These inequalities will allow us to normalize the auxiliary quantity $M$ in the proof of Lemma~\ref{lem:higher-pointwise} later on in Section~\ref{sec:technicallemmas}.

For $j\geq2$, write
\begin{equation*}
\Pi_j(\lambda)=\prod_{h=1}^{j-1}(1-h\lambda).
\end{equation*}
When $\lambda=1/m$, the quantity $\Pi_j(\lambda)$ is the labelled
$K_j$-density of the balanced complete $m$-partite graph.  Allowing
$\lambda$ to vary continuously gives the product curve that will be used
in the large-$u$ comparison in Lemma~\ref{lem:higher-product-comparisons}. 

\begin{lemma}\label{lem:higher-scalings}
For $0<x\leq 1$, both
\begin{equation}\label{eq:higher-scaling-quotients} \frac{\widehatFplus{\rtwotilde}{\rthreetilde}(x)} {x^{\rthreetilde/\rtwotilde}} \quad\text{and}\quad \frac{\widehatFplus{\rtwotilde}{\rthreetilde}(x)} {x^{(\rthreetilde-1)/(\rtwotilde-1)}}
\end{equation}
are nondecreasing, and they are strictly increasing on $\left[\Pi_{\rtwotilde}\left(\frac{1}{\rthreetilde-1}\right),1\right]$. Consequently, if $0<\lambda\leq1$, then
\begin{equation}\label{eq:higher-scaling}
\widehatFplus{\rtwotilde}{\rthreetilde} \bigl(\lambda^{\rtwotilde-1}x\bigr) \leq \lambda^{\rthreetilde-1} \widehatFplus{\rtwotilde}{\rthreetilde}(x).
\end{equation}
\end{lemma}

 \begin{proof} For $0<x<1$, by \eqref{eq:labelled-Fplus}, passing from $\Fplus{\rtwotilde}{\rthreetilde}$ to $\widehatFplus{\rtwotilde}{\rthreetilde}$ changes each quotient only by a positive constant. The two monotonicity assertions therefore follow from Lemma~\ref{lem:profile-monotonicities}. Applying the second one at $\lambda^{\rtwotilde-1}x$ and at $x$ gives \eqref{eq:higher-scaling}. We extend to $x=1$ by continuity.  \end{proof}
 
The next lemma is a higher-order form of Reiher's convex
product bound \cite[Proposition~3.1]{Reiher}.

\begin{lemma}[Higher-order convex product bound]
\label{lem:higher-product-bound}
If $0\leq\lambda\leq1/(\rthreetilde-1)$, then
\begin{equation}\label{eq:higher-product-bound}
\widehatFplus{\rtwotilde}{\rthreetilde}\bigl(\Pi_{\rtwotilde}(\lambda)\bigr)
\geq
\Pi_{\rthreetilde}(\lambda).
\end{equation}
\end{lemma}

\begin{proof}
If $\lambda=0$, then both sides of
\eqref{eq:higher-product-bound} are equal to $1$.  We may therefore
assume that
$$
0<\lambda\leq\frac1{\rthreetilde-1}.
$$
In particular, $\Pi_{\rtwotilde}(\lambda)<1$, so by
\eqref{eq:labelled-Fplus}, inequality
\eqref{eq:higher-product-bound} is equivalent to
$$
\Fplus{\rtwotilde}{\rthreetilde}
\left(
\frac{\Pi_{\rtwotilde}(\lambda)}{\rtwotilde!}
\right)
\geq
\frac{\Pi_{\rthreetilde}(\lambda)}{\rthreetilde!}.
$$
We prove this unlabelled form.

Recall from
Lemma~\ref{lem:profile-monotonicities} that the pair $ \bigl(D_{j_1}(N,\beta),D_{j_2}(N,\beta)\bigr) $ parametrizes the scallop between the balanced $N$- and $(N+1)$-partite points.  Direct differentiation gives

$$
\frac{d}{d\beta}D_j(N,\beta)
=
\binom{N}{j-1}
\frac{(j-1)(1-\beta)}{(N+\beta)^{j+1}}.
$$

Consequently, on a scallop where the $K_{\rthreetilde}$-density is positive, the chain rule yields

$$
\frac{dD_{\rthreetilde}}{dD_{\rtwotilde}}
=\frac{\frac{d}{d\beta}D_{\rthreetilde}(N,\beta)}{\frac{d}{d\beta}D_{\rtwotilde}(N,\beta)}=
\frac{(\rthreetilde-1)\binom{N}{\rthreetilde-1}}
{(\rtwotilde-1)\binom{N}{\rtwotilde-1}}
(N+\beta)^{\rtwotilde-\rthreetilde},
$$

which decreases with $\beta$.  By continuity at the endpoints,
every positive scallop of
$\Fplus{\rtwotilde}{\rthreetilde}$ is concave.

We next show that the product curve

$$
\frac{\Pi_{\rtwotilde}(\lambda)}{\rtwotilde!}
\longmapsto
\frac{\Pi_{\rthreetilde}(\lambda)}{\rthreetilde!}
$$

is convex.  It is enough to check the adjacent maps

$$
\frac{\Pi_j(\lambda)}{j!}
\longmapsto
\frac{\Pi_{j+1}(\lambda)}{(j+1)!}.
$$

Put

$$
S_j(\lambda)
=
\sum_{h=1}^{j-1}\frac{h}{1-h\lambda}.
$$

Since

$$
\frac{d\Pi_j(\lambda)}{d\lambda}
=-\Pi_j(\lambda)S_j(\lambda),
$$

we have again by the chain rule

\begin{align*}
\frac{d\bigl(\Pi_{j+1}(\lambda)/(j+1)!\bigr)}
{d\bigl(\Pi_j(\lambda)/j!\bigr)}
&=\frac{\frac{d\bigl(\Pi_{j+1}(\lambda)/(j+1)!\bigr)}{d\lambda}}{\frac{d\bigl(\Pi_{j}(\lambda)/(j!)\bigr)}{d\lambda}}=\frac{-\Pi_{j+1}(\lambda)S_{j+1}(\lambda)j!}{-\Pi_j(\lambda)S_j(\lambda)(j+1)!}\\
&=\frac{1-j\lambda}{j+1}\cdot \frac{S_{j+1}(\lambda)}{S_j(\lambda)}=\frac{(1-j\lambda)(S_j(\lambda)+\frac{j}{1-j\lambda})}{(j+1)S_j(\lambda)}\\
&=\frac{1}{j+1}
\left(
1-j\lambda+\frac{j}{S_j(\lambda)}
\right).
\end{align*}

This quantity is positive, and its derivative with respect to
$\lambda$ is

$$
-\frac{j}{j+1}
\left(
1+\frac{\frac{d S_j(\lambda)}{d\lambda}}{S_j(\lambda)^2}
\right)
<0
$$
because $d S_j(\lambda)/d\lambda=\sum_{h=1}^{j-1} h^2/(1-h\lambda)^2$.

Since $\frac{d\Pi_j(\lambda)}{d\lambda}<0$, the slope is increasing as a function of
$\Pi_j(\lambda)/j!$.  These properties extend to
$\lambda=1/(\rthreetilde-1)$ by continuity. Thus each adjacent map is increasing and convex,
and therefore so is their composition.

For every $N\geq\rthreetilde-1$, the profile scallop and the
product curve have the same two endpoints, because

$$
D_j(N,0)
=
\frac{\Pi_j(1/N)}{j!}
\qquad\text{and}\qquad
D_j(N,1)
=
\frac{\Pi_j(1/(N+1))}{j!}.
$$

The concave profile scallop lies above the chord joining these endpoints,
while the convex product curve lies below it.  Hence

$$
\Fplus{\rtwotilde}{\rthreetilde}
\left(\frac{\Pi_{\rtwotilde}(\lambda)}{\rtwotilde!}\right)
\geq
\frac{\Pi_{\rthreetilde}(\lambda)}{\rthreetilde!}
$$
whenever
$$
\frac1{N+1}\leq\lambda\leq\frac1N,
\qquad
N\geq\rthreetilde-1.
$$
These intervals cover
$(0,1/(\rthreetilde-1)]$.  Together with the already verified case
$\lambda=0$, this proves \eqref{eq:higher-product-bound}.
\end{proof}

\begin{remark}
When $p=3$, we have $\rtwotilde=2$ and
$\Pi_2(\lambda)=1-\lambda$.  Writing $z=1-\lambda$, the conclusion becomes
$$
\widehatFplus{2}{\rthreetilde}(z)
\geq z(2z-1)(3z-2)\cdots((\rthreetilde-1)z-(\rthreetilde-2)),
$$
which is exactly Reiher's product inequality
\cite[Proposition~3.1]{Reiher}.
\end{remark}

\section{The two link comparisons}
\label{sec:link-comparisons}

The pointwise link lower bound for $\frac{\rtwotilde KD-P}{\rtwotilde-1}$ in Lemma~\ref{lem:higher-pointwise} divides naturally at
$$
u=\frac{\rthreetilde}{\rthreetilde-\rtwotilde}.
$$
The small-$u$ range is a comparison with the multipartite
$\{K_{\rtwotilde},K_{\rthreetilde}\}$-profile.  The large-$u$ range is
controlled by the higher-order product inequality.  These are the two
comparisons needed in the pointwise proof.

Throughout this section, put
\begin{equation}\label{eq:higher-Y}
Y=(s-\rtwotilde)_{\rthreetilde-\rtwotilde},
\qquad
K=(s-1)_{\rtwotilde-1}.
\end{equation}

\begin{lemma}[Higher-order small-$u$ comparison]
\label{lem:higher-small-u}
Let $s\geq\rthreetilde$, let $Y,K$ be given by
\eqref{eq:higher-Y}, and suppose
$$
1\leq u\leq\frac{\rthreetilde}{\rthreetilde-\rtwotilde}.
$$
Put
\begin{equation}\label{eq:higher-y0}
y_0(u)=\frac{s}{\rtwotilde}
\left(\rtwotilde-1+\frac{\rthreetilde-\rtwotilde}{\rthreetilde}u\right).
\end{equation}
Then
\begin{equation}\label{eq:higher-small-u}
\widehatFplus{\rtwotilde}{\rthreetilde}\left(
\frac{(\rthreetilde-\rtwotilde)sKu}{\rthreetilde\,y_0(u)^{\rtwotilde}}
\right)
\geq
\frac{Y(\rthreetilde-\rtwotilde)sK(u-1)}
{\rtwotilde\,y_0(u)^{\rthreetilde}}.
\end{equation}
\end{lemma}

\begin{proof}
Put
$$
w=\frac{1-\frac{\rthreetilde-\rtwotilde}{\rthreetilde}u}
{\rtwotilde-1+\frac{\rthreetilde-\rtwotilde}{\rthreetilde}u}.
$$
To relate these expressions to the notation
$D_j(N,\beta)$ introduced in
Lemma~\ref{lem:profile-monotonicities}, put
$$
N=s-1,
\qquad
\beta=\frac{1-(s-1)w}{1+w}.
$$
Recall that when we introduced $D_j(N,\beta)$, we defined it for $0\leq \beta \leq 1$. Here, we allow $\beta$ to be negative.  Since
$$
N+\beta=\frac{s}{1+w},
$$
a direct calculation gives, for every $j\geq2$,
\begin{equation}\label{eq:higher-Dj-bridge}
j!D_j(s-1,\beta)
=
\Pi_j(1/s)
\bigl(1-(j-1)w\bigr)(1+w)^{j-1}.
\end{equation}

Note that $0\leq w\leq1/(\rthreetilde-1)$, and direct substitution in
\eqref{eq:higher-Y}--\eqref{eq:higher-y0} gives

$$\frac{(\rthreetilde-\rtwotilde)sKu}{\rthreetilde\,y_0(u)^{\rtwotilde}}=\Pi_{\rtwotilde}(1/s)(1-(\rtwotilde-1)w)(1+w)^{\rtwotilde-1},$$
and
$$\frac{Y(\rthreetilde-\rtwotilde)sK(u-1)}{\rtwotilde\,y_0(u)^{\rthreetilde}}=\Pi_{\rthreetilde}(1/s)(1-(\rthreetilde-1)w)(1+w)^{\rthreetilde-1}.$$
By \eqref{eq:higher-Dj-bridge}, the two right-hand sides are
respectively
$$
\rtwotilde!D_{\rtwotilde}(s-1,\beta)
\qquad\text{and}\qquad
\rthreetilde!D_{\rthreetilde}(s-1,\beta).
$$
If $w\leq1/(s-1)$, then $0\leq\beta\leq1$, so these are the labelled
densities of a genuine complete multipartite profile point.  If
$w>1/(s-1)$, then
$$
-\frac{s-\rthreetilde}{\rthreetilde}
\leq\beta<0,
$$
and the same formulas describe the formal negative-exceptional-part
continuation covered by Li--Liu--Zhang's compression Lemma
\ref{lem:compression}.  In either case,
$$
\widehatFplus{\rtwotilde}{\rthreetilde}
\bigl(\rtwotilde!D_{\rtwotilde}(s-1,\beta)\bigr)
\geq
\rthreetilde!D_{\rthreetilde}(s-1,\beta),
$$
which is exactly \eqref{eq:higher-small-u}.
\end{proof}

The following product comparison supplies the large-$u$ interpolation.

\begin{lemma}[Higher-order product-shift comparison]
\label{lem:higher-product-comparisons}
Suppose $s\geq\rthreetilde$ and
$$
0\leq x\leq\frac1{(\rtwotilde-1)(\rthreetilde-1)}.
$$
If $\lambda_x\in \left[\frac{1}{s},\frac{1}{\rtwotilde-1}\right]$ is determined by
\begin{equation*}
\Pi_{\rtwotilde}(\lambda_x)
=
\Pi_{\rtwotilde}(1/s)
(1+(\rtwotilde-1)x)(1-x)^{\rtwotilde-1},
\end{equation*}
then
\begin{equation}\label{eq:higher-product-shift}
\begin{aligned}
\frac{\Pi_{\rthreetilde}(\lambda_x)}{\Pi_{\rthreetilde}(1/s)}
\geq{}&
(1+(\rthreetilde-1)x)(1-x)^{\rthreetilde-1}\\
&\times
\bigl(1-(\rtwotilde-1)(\rthreetilde-1)x\bigr)^{(\rthreetilde-\rtwotilde)/(\rtwotilde-1)}.
\end{aligned}
\end{equation}
\end{lemma}

\begin{proof}
Write
$$
g_j(x)=(1+(j-1)x)(1-x)^{j-1},
\qquad
\mu=\frac1s+\left(1-\frac1s\right)x^{2},
$$
and note that both $g_j(x)$ and $\mu$ are nonnegative in our setup. Then $\lambda_x$ is determined by
$\Pi_{\rtwotilde}(\lambda_x)=\Pi_{\rtwotilde}(1/s)g_{\rtwotilde}(x)$ and the
conclusion \eqref{eq:higher-product-shift} reads
$$
\frac{\Pi_{\rthreetilde}(\lambda_x)}{\Pi_{\rthreetilde}(1/s)}
\geq
g_{\rthreetilde}(x)
\bigl(1-(\rtwotilde-1)(\rthreetilde-1)x\bigr)^{(\rthreetilde-\rtwotilde)/(\rtwotilde-1)}.
$$
We will use throughout that
\begin{equation}\label{eq:g-log-derivative}
\frac{g_j'(x)}{g_j(x)}
=-\frac{j(j-1)x}{(1+(j-1)x)(1-x)},
\end{equation}
which is immediate from
$g_j'(x)=-j(j-1)x(1-x)^{j-2}$ which is nonpositive, and negative for $x>0$.

\medskip
\emph{Step 0:} We first show that $\lambda_x$ is well defined.
For $0\leq\lambda\leq1/(\rtwotilde-1)$ every factor of
$\Pi_{\rtwotilde}(\lambda)$ is nonnegative and strictly decreasing in
$\lambda$, so $\Pi_{\rtwotilde}$ decreases strictly from $1$ to $0$ on this
interval.  We have that $g_{\rtwotilde}(0)=1$, and by \eqref{eq:g-log-derivative}, we have $g_{\rtwotilde}'\leq0$, so $0\leq g_{\rtwotilde}(x)\leq1$.  Hence
$\Pi_{\rtwotilde}(1/s)g_{\rtwotilde}(x)$ lies in
$[0,\Pi_{\rtwotilde}(1/s)]$, and there is a unique
$$
\lambda_x\in\left[\frac1s,\frac1{\rtwotilde-1}\right]
$$
with $\Pi_{\rtwotilde}(\lambda_x)=\Pi_{\rtwotilde}(1/s)g_{\rtwotilde}(x)$.

\medskip
\emph{Step 1: We show that $\mu\leq1/(\rthreetilde-1)$.}
Since $x^{2}<1$, the quantity $\mu$ is decreasing in
$s$. Using $s\geq\rthreetilde$ and $x\leq \frac1{(\rtwotilde-1)(\rthreetilde-1)}$ gives
$$
\mu\leq
\frac1\rthreetilde+\left(1-\frac1\rthreetilde\right)x^{2}
\leq
\frac1\rthreetilde+
\left(1-\frac1\rthreetilde\right)
\frac1{(\rtwotilde-1)^{2}(\rthreetilde-1)^{2}}
=\frac1\rthreetilde+
\frac1{\rthreetilde(\rtwotilde-1)^{2}(\rthreetilde-1)}\leq \frac1{\rthreetilde-1},
$$
where the last inequality follows from $\rtwotilde\geq2$.

\medskip
\emph{Step 2:} We now show that
$\lambda_x-\frac1s\leq\bigl(1-\frac1s\bigr)x^{2}$.
This is equivalent to showing that $\lambda_x\leq\mu$. We know that $0\leq \mu \leq \frac{1}{\rthreetilde-1}\leq \frac{1}{\rtwotilde-1}$ and that $\frac{1}{s} \leq \lambda_x \leq \frac{1}{\rtwotilde -1}$ and that
$\Pi_{\rtwotilde}$ is strictly decreasing on $[0,1/(\rtwotilde-1)]$.  It therefore suffices to prove that
$$
\Pi_{\rtwotilde}(\mu)\leq\Pi_{\rtwotilde}(\lambda_x)
=\Pi_{\rtwotilde}(1/s)g_{\rtwotilde}(x).
$$
For $1\leq h\leq\rtwotilde-1$, we have $h/s<1$ and hence
$\frac{1-1/s}{1-h/s}\geq1$, so
$$
1-h\mu
=\left(1-\frac hs\right)
\left(1-\frac{h(1-1/s)x^{2}}{1-h/s}\right)
\leq\left(1-\frac hs\right)(1-hx^{2}).
$$
All factors here are nonnegative, because
$hx^{2}\leq(\rtwotilde-1)x^{2}\leq\frac1{(\rtwotilde-1)(\rthreetilde-1)^{2}}<1$.
Taking the product over $1\leq h\leq\rtwotilde-1$ on both sides of the previous line gives
$$
\Pi_{\rtwotilde}(\mu)
\leq
\Pi_{\rtwotilde}(1/s)\prod_{h=1}^{\rtwotilde-1}(1-hx^{2}),
$$
so it is enough to show
\begin{equation}\label{eq:square-product}
\prod_{h=1}^{\rtwotilde-1}(1-hx^{2})
\leq g_{\rtwotilde}(x).
\end{equation}
If $\rtwotilde=2$, both sides equal $1-x^{2}$ and
\eqref{eq:square-product} thus holds.  Suppose $\rtwotilde\geq3$ and take logs on both sides of \eqref{eq:square-product}. We will show that $$
R(x):=\log g_{\rtwotilde}(x)-\sum_{h=1}^{\rtwotilde-1}\log(1-hx^{2})\geq 0.
$$
  By \eqref{eq:g-log-derivative},
$$
R'(x)=
\sum_{h=1}^{\rtwotilde-1}\frac{2hx}{1-hx^{2}}
-\frac{\rtwotilde(\rtwotilde-1)x}
{(1+(\rtwotilde-1)x)(1-x)}.
$$
Since $0<1-hx^{2}\leq1$, the sum is at least
$\sum_{h=1}^{\rtwotilde-1}2hx=\rtwotilde(\rtwotilde-1)x$.  Now note that
$$
(1+(\rtwotilde-1)x)(1-x)-1
=x\bigl((\rtwotilde-2)-(\rtwotilde-1)x\bigr)
$$
by simply expanding. Furthermore, 
$\rtwotilde\geq3$ and
$x\leq\frac1{(\rtwotilde-1)(\rthreetilde-1)}
\leq\frac{\rtwotilde-2}{\rtwotilde-1}$, and so the right-hand side of the previous equation is nonnegative, which implies that the denominator of the subtracted term in the expression for $R'(x)$ is at least $1$.  Hence that subtracted term is at
most $\rtwotilde(\rtwotilde-1)x$, so $R'(x)\geq0$ on the stated range and
\eqref{eq:square-product} follows.  This proves
\begin{equation}\label{eq:higher-lambda-shift}
\lambda_x-\frac1s\leq\left(1-\frac1s\right)x^2.
\end{equation}
Moreover, with the fact that we showed that $\mu\leq 1/(\rthreetilde-1)$ in Step 1, this also
records the bound

$$\lambda_x\leq\frac1{\rthreetilde-1}$$
that will be needed in Section~\ref{sec:technicallemmas}.

\medskip
\emph{Step 3:}
We claim that
\begin{equation}\label{eq:tail-coefficient}
\frac{h(1-1/s)}{1-h/s}=\frac{h(s-1)}{s-h}\leq(\rthreetilde-1)^{2}
\qquad\text{for }\rtwotilde\leq h\leq\rthreetilde-1.
\end{equation}
The left side increases in $h$, so it suffices to treat
$h=\rthreetilde-1$, where the claim reads
$(s-1)\leq(\rthreetilde-1)(s-\rthreetilde+1)$.  Rearranging, this is
$$
(\rthreetilde-2)(\rthreetilde)\leq(\rthreetilde-2)s ,
$$
which holds because $s\geq\rthreetilde$.

\medskip
\emph{Step 4: } We give a first lower bound for $\frac{\Pi_{\rthreetilde}(\lambda_x)}{\Pi_{\rthreetilde}(1/s)}$.
Split the quotient at $h=\rtwotilde-1$:
$$
\frac{\Pi_{\rthreetilde}(\lambda_x)}{\Pi_{\rthreetilde}(1/s)}
=
\underbrace{\frac{\Pi_{\rtwotilde}(\lambda_x)}{\Pi_{\rtwotilde}(1/s)}}
_{=\,g_{\rtwotilde}(x)}
\cdot
\prod_{h=\rtwotilde}^{\rthreetilde-1}\frac{1-h\lambda_x}{1-h/s},
$$
the first factor being $g_{\rtwotilde}(x)$ by the definition of
$\lambda_x$.  For each remaining index, 
$$
\frac{1-h\lambda_x}{1-h/s}
=1-\frac{h(\lambda_x-1/s)}{1-h/s}
\geq1-\frac{h(1-1/s)x^{2}}{1-h/s}
\geq1-(\rthreetilde-1)^{2}x^{2}\geq0,
$$
where the first inequality follows from Step~2, the second inequality follows from
\eqref{eq:tail-coefficient}, and the last inequality follows from
$(\rthreetilde-1)^{2}x^{2}\leq(\rtwotilde-1)^{-2}\leq1$.  There are
$\rthreetilde-\rtwotilde$ such indices, so
\begin{equation}\label{eq:higher-product-shift-rough}
\frac{\Pi_{\rthreetilde}(\lambda_x)}{\Pi_{\rthreetilde}(1/s)}
\geq
g_{\rtwotilde}(x)
\bigl(1-(\rthreetilde-1)^{2}x^{2}\bigr)^{\rthreetilde-\rtwotilde}.
\end{equation}

\medskip
\emph{Step 5:} Finally, we will show how \eqref{eq:higher-product-shift-rough} implies
\eqref{eq:higher-product-shift}. We claim that it suffices to
prove
\begin{equation}\label{eq:step5mid}
g_{\rtwotilde}(x)
\bigl(1-(\rthreetilde-1)^{2}x^{2}\bigr)^{\rthreetilde-\rtwotilde}
\geq
g_{\rthreetilde}(x)
\bigl(1-(\rthreetilde-1)x\bigr)^{\rthreetilde-\rtwotilde}.
\end{equation}

Indeed, in that case, we would obtain the desired conclusion of our lemma that
\begin{align*}
\frac{\Pi_{\rthreetilde}(\lambda_x)}{\Pi_{\rthreetilde}(1/s)}
\geq
g_{\rtwotilde}(x)
\bigl(1-(\rthreetilde-1)^{2}x^{2}\bigr)^{\rthreetilde-\rtwotilde}&\geq
g_{\rthreetilde}(x)
\bigl(1-(\rthreetilde-1)x\bigr)^{\rthreetilde-\rtwotilde}\\
&\geq g_{\rthreetilde}(x)
\bigl(1-(\rtwotilde-1)(\rthreetilde-1)x\bigr)^{(\rthreetilde-\rtwotilde)/(\rtwotilde-1)}
\end{align*}
where the first inequality follows from \eqref{eq:higher-product-shift-rough}, the second inequality is our claim, and the last inequality is an application of Bernoulli's inequality. Indeed, since $(\rthreetilde-1)x\leq\frac1{\rtwotilde-1}\leq1$, Bernoulli's
inequality applied with exponent $\rtwotilde-1$ gives
$$
1-(\rtwotilde-1)(\rthreetilde-1)x
\leq\bigl(1-(\rthreetilde-1)x\bigr)^{\rtwotilde-1},
$$
and raising both sides to the positive power
$(\rthreetilde-\rtwotilde)/(\rtwotilde-1)$ yields
$$
\bigl(1-(\rtwotilde-1)(\rthreetilde-1)x\bigr)
^{(\rthreetilde-\rtwotilde)/(\rtwotilde-1)}
\leq
\bigl(1-(\rthreetilde-1)x\bigr)^{\rthreetilde-\rtwotilde}.
$$

So we now want to prove \eqref{eq:step5mid}. Factoring
$1-(\rthreetilde-1)^{2}x^{2}
=(1-(\rthreetilde-1)x)(1+(\rthreetilde-1)x)$ and cancelling
$(1-(\rthreetilde-1)x)^{\rthreetilde-\rtwotilde}$, which is legitimate
since that quantity is nonnegative and the inequality is trivial when it
vanishes, we are left with
$$
(1+(\rtwotilde-1)x)(1-x)^{\rtwotilde-1}
(1+(\rthreetilde-1)x)^{\rthreetilde-\rtwotilde}
\geq
(1+(\rthreetilde-1)x)(1-x)^{\rthreetilde-1}.
$$
Dividing by the positive quantities $(1-x)^{\rtwotilde-1}$ and
$1+(\rthreetilde-1)x$ turns this into
$$
(1+(\rtwotilde-1)x)
(1+(\rthreetilde-1)x)^{\rthreetilde-\rtwotilde-1}
\geq(1-x)^{\rthreetilde-\rtwotilde}.
$$
Here $\rthreetilde-\rtwotilde-1\geq0$, so the left side is at least $1$,
while the right side is at most $1$.  This completes the proof.
\end{proof}

\section{Proof of the two technical lemmas}
\label{sec:technicallemmas}

We are finally ready to prove the two technical lemmas. Note that Lemmas~3.3 and~3.4 in \cite{LLZ} prove the small-$u$ part of the next
lemma; the new higher-order product bound and
Lemma~\ref{lem:higher-product-comparisons} allow us to prove the large-$u$ part.

\begin{proof}[Proof of Lemma~\ref{lem:higher-pointwise}]
Set
$$
L=\frac{\rtwotilde KD-P}{\rtwotilde-1},
$$
i.e., the functional described in Section~\ref{sec:generalideatechnicalstatements},
and let
$$
X=(\rthreetilde-1)P-(\rthreetilde-\rtwotilde)KD-c
$$
be the factor appearing on the right side of the KKT vertex expression
\eqref{eq:higher-local-profile}.  We want a variable $M$ that normalizes
the candidate supporting line and a variable $u$ for which $u-1$
normalizes $X$.  Suppose that
$$
sKM=L+\beta c,
\qquad
\eta=\frac{\kappa P}{sK},
\qquad
u=\frac\eta M,
$$
and require
$$
sKM(u-1)=\chi X
$$
for constants $\beta,\kappa,\chi$.  Expanding both sides and comparing
the coefficients of $D$, $P$, and $c$ forces
$$
\chi=\beta
=\frac{\rtwotilde}
{(\rtwotilde-1)(\rthreetilde-\rtwotilde)},
\qquad
\kappa=\frac{\rthreetilde}{\rthreetilde-\rtwotilde}.
$$
Thus the normalizations below are forced by the affine functional and the
KKT slack, rather than guessed.  Put
$$
M=\frac{L+\widehat c}{sK},
\qquad
\eta=\frac{\rthreetilde P}{(\rthreetilde-\rtwotilde)sK},
\qquad
u=\frac{\eta}{M},
$$
and
\begin{equation*}
y=\frac{D}{M^{1/\rtwotilde}},
\qquad
y_0:=y_0(u)=\frac{s}{\rtwotilde}
\left(\rtwotilde-1+\frac{\rthreetilde-\rtwotilde}{\rthreetilde}u\right).
\end{equation*}
We have that $M>0$ if and only if $P<\rtwotilde K D+\frac{\rtwotilde}{\rthreetilde-\rtwotilde}c$ which we know is true by the KKT edge inequality~\eqref{eq:higher-local-edge} $P\leq K D+\frac{\rthreetilde}{(\rthreetilde-1)(\rthreetilde-\rtwotilde)}c<\rtwotilde KD+\frac{\rtwotilde}{\rthreetilde-\rtwotilde}c$, where the strict inequality follows from the fact that $D>0$, $\rtwotilde\geq 2$ and $\frac{\rthreetilde}{\rthreetilde-1}<2\leq \rtwotilde$.   Moreover, $u\geq 1$ since
\begin{equation*}
\eta-M
=
\frac{\rtwotilde((\rthreetilde-1)P-(\rthreetilde-\rtwotilde)KD-c)}
{(\rtwotilde-1)(\rthreetilde-\rtwotilde)sK}
\geq0
\end{equation*}
because the right side of the KKT vertex inequality \eqref{eq:higher-local-profile} is nonnegative and thus $(\rthreetilde-1)P-(\rthreetilde-\rtwotilde)KD-c\geq 0$.  Using again this last inequality, we have that \eqref{eq:higher-local-profile} implies
\begin{equation}\label{eq:higher-normalized-local}
\widehatFplus{\rtwotilde}{\rthreetilde}\left(
\frac{(\rthreetilde-\rtwotilde)sKu}{\rthreetilde\,y^{\rtwotilde}}
\right)
\leq
\frac{Y(\rthreetilde-\rtwotilde)sK(u-1)}
{\rtwotilde M^{(\rthreetilde-\rtwotilde)/\rtwotilde}y^{\rthreetilde}}.
\end{equation}
since $\frac{P}{D^{\rtwotilde}}=\frac{(\rthreetilde-\rtwotilde)sKu}{\rthreetilde\,y^{\rtwotilde}}$.  The definitions also give
\begin{equation}\label{eq:higher-central-identity}
D+\frac{(\rtwotilde-1)\widehat c}{\rtwotilde K}=My_0.
\end{equation}

The desired pointwise inequality is equivalent to
$$
M^{(\rtwotilde-1)/\rtwotilde}
\geq
1+\frac{(\rtwotilde-1)\widehat c}{\rtwotilde KD}.
$$
Using~\eqref{eq:higher-central-identity}, this is in turn equivalent to
$$
y\geq y_0.
$$

For brevity, put
$$
Z=
\frac{(\rthreetilde-\rtwotilde)sKu}{\rthreetilde},
\qquad
W=
\frac{Y(\rthreetilde-\rtwotilde)sK(u-1)}{\rtwotilde}.
$$
Thus \eqref{eq:higher-normalized-local} becomes
$$
\widehatFplus{\rtwotilde}{\rthreetilde}
\left(\frac{Z}{y^{\rtwotilde}}\right)
\leq
\frac{W}
{M^{(\rthreetilde-\rtwotilde)/\rtwotilde}
 y^{\rthreetilde}}.
$$

We prove that $y\geq y_0$ in two ranges.  The dividing value
$$
u=\frac{\rthreetilde}{\rthreetilde-\rtwotilde}
$$
is the value of $u$ for the exceptional link of the candidate graph when
$c=0$ and $M=1$.  Below it, the critical pair in
Lemma~\ref{lem:higher-small-u} can be parametrized by the multipartite
$\{K_{\rtwotilde},K_{\rthreetilde}\}$-profile.  Depending on the parameters, this is either a genuine profile
point or its formal negative-exceptional-part continuation.
In the latter case, Li--Liu--Zhang's compression Lemma~\ref{lem:compression} still places it below the genuine profile. Above the dividing value, that parametrization
leaves the useful range, and the higher-order product inequality replaces
it.

Suppose first that
$$
1<u\leq
\frac{\rthreetilde}{\rthreetilde-\rtwotilde}.
$$
We first show that $M\geq1$.  Suppose instead that $M<1$, and put
$$
\widetilde D=\frac DM,
\qquad
\widetilde P=\frac PM,
\qquad
\widetilde c=\frac cM.
$$
Thus $0\leq \widetilde P \leq \widetilde D^{\rtwotilde}$. If $z=P/D^{\rtwotilde}$, then
\eqref{eq:higher-scaling} gives
$$
\begin{aligned}
\widetilde D^{\rthreetilde}
\widehatFplus{\rtwotilde}{\rthreetilde}
\left(
\frac{\widetilde P}{\widetilde D^{\rtwotilde}}
\right)
&=
\frac{D^{\rthreetilde}}{M^{\rthreetilde}}
\widehatFplus{\rtwotilde}{\rthreetilde}
\left(M^{\rtwotilde-1}z\right)
\\
&\leq
\frac1M D^{\rthreetilde}
\widehatFplus{\rtwotilde}{\rthreetilde}(z).
\end{aligned}
$$
The right side of \eqref{eq:higher-local-profile} is also divided by
$M$, and \eqref{eq:higher-local-edge} scales in exactly the same way.
Thus $(\widetilde D,\widetilde P,\widetilde c)$ still satisfies the
two local inequalities.

Put
$$
\widehat{\widetilde c}
=
\frac{\rtwotilde\widetilde c}
{(\rtwotilde-1)(\rthreetilde-\rtwotilde)}.
$$
The corresponding normalized parameters satisfy
$$
\widetilde M=1,
\qquad
\widetilde u=u.
$$
Dividing \eqref{eq:higher-central-identity} by $M$ gives
$$
\widetilde D+
\frac{(\rtwotilde-1)\widehat{\widetilde c}}
{\rtwotilde K}
=
y_0,
$$
and hence $\widetilde D<y_0$.

On the other hand, the normalized local inequality for the scaled
triple and Lemma~\ref{lem:higher-small-u} give
$$
\frac{
\widehatFplus{\rtwotilde}{\rthreetilde}
(Z/\widetilde D^{\rtwotilde})
}{
(Z/\widetilde D^{\rtwotilde})^{
\rthreetilde/\rtwotilde}
}
\leq
\frac{W}{Z^{\rthreetilde/\rtwotilde}}
\leq
\frac{
\widehatFplus{\rtwotilde}{\rthreetilde}
(Z/y_0^{\rtwotilde})
}{
(Z/y_0^{\rtwotilde})^{
\rthreetilde/\rtwotilde}
}.
$$
Since $\widetilde D<y_0$, the argument on the left is strictly larger
than the argument on the right.  This contradicts the strict
monotonicity of the first quotient in
\eqref{eq:higher-scaling-quotients} (monotonicity is strict since $u>1$, which gives that $W>0$, and Lemma~\ref{lem:higher-small-u} gives positivity at the comparison point).  Therefore $M\geq1$.

If $y<y_0$, then the same comparison gives
$$
\begin{aligned}
\frac{
\widehatFplus{\rtwotilde}{\rthreetilde}
(Z/y^{\rtwotilde})
}{
(Z/y^{\rtwotilde})^{\rthreetilde/\rtwotilde}
}
&\leq
\frac{W}
{M^{(\rthreetilde-\rtwotilde)/\rtwotilde}
 Z^{\rthreetilde/\rtwotilde}}
\\
&\leq
\frac{W}{Z^{\rthreetilde/\rtwotilde}}
\\
&\leq
\frac{
\widehatFplus{\rtwotilde}{\rthreetilde}
(Z/y_0^{\rtwotilde})
}{
(Z/y_0^{\rtwotilde})^{\rthreetilde/\rtwotilde}
},
\end{aligned}
$$
again contradicting strict monotonicity.  Hence
\begin{equation}\label{eq:higher-y-ge-y0}
y\geq y_0.
\end{equation}

When $u=1$, we can use Li--Liu--Zhang's Lemmas~3.3-3.4 \cite{LLZ}. We include the details here for completeness. When $u=1$, then $W=0$, so
\eqref{eq:higher-normalized-local} forces
$$
\widehatFplus{\rtwotilde}{\rthreetilde}
\left(\frac{Z}{y^{\rtwotilde}}\right)
\leq0,
$$
and hence $\widehatFplus{\rtwotilde}{\rthreetilde}(Z/y^{\rtwotilde})=0$.
We now show that
$$
\frac{Z}{y^{\rtwotilde}}
\leq
\Pi_{\rtwotilde}\!\left(\frac1{\rthreetilde-1}\right).
$$

By Lemma~\ref{lem:profile-monotonicities}, the function
$x\mapsto\Fplus{\rtwotilde}{\rthreetilde}(x)/x^{\rthreetilde/\rtwotilde}$
is nondecreasing, and hence so is
$\Fplus{\rtwotilde}{\rthreetilde}$ itself; by
\eqref{eq:labelled-Fplus} the same holds for
$\widehatFplus{\rtwotilde}{\rthreetilde}$.  Therefore $\widehatFplus{\rtwotilde}{\rthreetilde}(x)=0$ on an initial interval, and it suffices to locate its right
endpoint.

On the scallop indexed by $s'\geq\rtwotilde-1$, the function
$\Fplus{\rtwotilde}{\rthreetilde}$ is given by
$$
\frac{1}{(s'+1)^{\rthreetilde}}
\binom{s'+1}{\rthreetilde}
(1+\alpha)^{\rthreetilde-1}
\bigl(1-(\rthreetilde-1)\alpha\bigr),
\qquad 0\leq\alpha\leq\frac1{s'} .
$$
If $s'\leq\rthreetilde-2$, the binomial coefficient vanishes and the
whole scallop is zero.  If $s'=\rthreetilde-1$, the expression vanishes
only at the left endpoint $\alpha=1/s'=1/(\rthreetilde-1)$ and is
positive for $\alpha<1/(\rthreetilde-1)$.  That endpoint is the balanced
complete $(\rthreetilde-1)$-partite graph, whose labelled
$K_{\rtwotilde}$-density is
$$
\rtwotilde!\,
\frac{\binom{\rthreetilde-1}{\rtwotilde}}
{(\rthreetilde-1)^{\rtwotilde}}
=
\frac{(\rthreetilde-2)_{\rtwotilde-1}}{(\rthreetilde-1)^{\rtwotilde-1}}
=
\Pi_{\rtwotilde}\!\left(\frac1{\rthreetilde-1}\right).
$$
Note that this is why this $u=1$ is the right threshold:  the right-hand side of \eqref{eq:higher-local-profile} vanishes, so the vertex link graph is forced to be $K_{\rthreetilde}$-free, and the extremal $K_{\rthreetilde}$-free vertex link graph is the balanced complete $(\rthreetilde-1)$-partite graph.

So since $\widehatFplus{\rtwotilde}{\rthreetilde}(Z/y^{\rtwotilde})=0$, 

$$\frac{Z}{y^{\rtwotilde}}
\leq
\Pi_{\rtwotilde}\!\left(\frac1{\rthreetilde-1}\right),$$ using $K=s^{\rtwotilde-1}\Pi_{\rtwotilde}(1/s)$ and $
Z=\frac{(\rthreetilde-\rtwotilde)sK}{\rthreetilde}$ at $u=1$, we have \begin{equation}\label{eq:yetanothermiddlestep}y^{\rtwotilde}
\geq
\frac{\rthreetilde-\rtwotilde}{\rthreetilde}\cdot s^{\rtwotilde} \cdot \frac{\Pi_{\rtwotilde}\left(\frac1s\right)}{\Pi_{\rtwotilde}\left(\frac1{\rthreetilde-1}\right)}.
\end{equation}
At $u=1$, we also have
$$
y_0(1)
=\frac{s}{\rtwotilde}
\left(\rtwotilde-1+\frac{\rthreetilde-\rtwotilde}{\rthreetilde}\right)
=\frac{s(\rthreetilde-1)}{\rthreetilde}.
$$
If we can show that 
\begin{equation}\label{eq:u1-reduction}
\frac{\Pi_{\rtwotilde}(1/s)}
{\Pi_{\rtwotilde}\bigl(1/(\rthreetilde-1)\bigr)}
\geq
\frac{(\rthreetilde-1)^{\rtwotilde}}
{\rthreetilde^{\,\rtwotilde-1}(\rthreetilde-\rtwotilde)}, 
\end{equation}
then we will have that the right-hand side of~\eqref{eq:yetanothermiddlestep} is greater or equal to $y_0(1)^{\rtwotilde}$, and so we will have $y\geq y_0(1)$ as desired.
Since $s\geq\rthreetilde$ and $\Pi_{\rtwotilde}$ is decreasing, we have
$\Pi_{\rtwotilde}(1/s)\geq\Pi_{\rtwotilde}(1/\rthreetilde)$, and
$$
\frac{\Pi_{\rtwotilde}(1/\rthreetilde)}
{\Pi_{\rtwotilde}\bigl(1/(\rthreetilde-1)\bigr)}
=
\frac{(\rthreetilde-1)_{\rtwotilde-1}}
{(\rthreetilde-2)_{\rtwotilde-1}}
\cdot
\frac{(\rthreetilde-1)^{\rtwotilde-1}}
{\rthreetilde^{\,\rtwotilde-1}}
=
\frac{\rthreetilde-1}{\rthreetilde-\rtwotilde}
\cdot
\frac{(\rthreetilde-1)^{\rtwotilde-1}}
{\rthreetilde^{\,\rtwotilde-1}}
=
\frac{(\rthreetilde-1)^{\rtwotilde}}
{\rthreetilde^{\,\rtwotilde-1}(\rthreetilde-\rtwotilde)},
$$
which is exactly the right-hand side of \eqref{eq:u1-reduction}.  Thus
\eqref{eq:u1-reduction} holds, with equality precisely when
$s=\rthreetilde$, and \eqref{eq:higher-y-ge-y0} holds at $u=1$ as well.

It remains to consider $u>\rthreetilde/(\rthreetilde-\rtwotilde)$.  The edge inequality gives
\begin{equation}\label{eq:higher-E}
\frac{y}{M^{(\rtwotilde-1)/\rtwotilde}}
\leq
E(u):=
\frac{s}{\rtwotilde\rthreetilde}
\left(
\rthreetilde^2(\rtwotilde-1)
-(\rthreetilde-\rtwotilde)(\rthreetilde(\rtwotilde-1)-\rtwotilde)u
\right).
\end{equation}
In particular,
$$
u<
\frac{\rthreetilde^2(\rtwotilde-1)}
{(\rthreetilde-\rtwotilde)(\rthreetilde(\rtwotilde-1)-\rtwotilde)}.
$$
Put
\begin{equation}\label{eq:vxdefs}
v=\frac{\rthreetilde-\rtwotilde}{\rthreetilde}u,
\qquad
x=\frac{v-1}{v+\rtwotilde-1}.
\end{equation}
Then substituting these into the
definitions of $Z$, $W$, $E$, and $y_0$ yields
\begin{equation}\label{eq:higher-large-x}
Z=\frac{sK\bigl(1+(\rtwotilde-1)x\bigr)}{1-x},
W=\frac{YsK\bigl(1+(\rthreetilde-1)x\bigr)}{1-x},
E(u)=\frac{s\bigl(1-(\rtwotilde-1)(\rthreetilde-1)x\bigr)}{1-x},
y_0(u)=\frac{s}{1-x}
\end{equation}
and we also have that  $0<x<\frac1{(\rtwotilde-1)(\rthreetilde-1)}$.

Suppose, for a contradiction, that $y<y_0$.  

First note that since $s\geq\rthreetilde$, the definitions of $Y$ and $K$ give
\begin{equation}\label{eq:higher-PiYK}
\Pi_{\rtwotilde}(1/s)
=\frac{K}{s^{\rtwotilde-1}},
\qquad
\Pi_{\rthreetilde}(1/s)
=\frac{YK}{s^{\rthreetilde-1}}.
\end{equation}
Define $\lambda_x$ by
$$
\Pi_{\rtwotilde}(\lambda_x)
=\frac{Z}{y_0^{\rtwotilde}}
=\Pi_{\rtwotilde}(1/s)
(1+(\rtwotilde-1)x)(1-x)^{\rtwotilde-1}.
$$

Moreover,
$$
\frac{Z}{y^{\rtwotilde}}
=
\frac{P}{D^{\rtwotilde}}
\leq1.
$$
Since $y<y_0$, we have
$$
\Pi_{\rtwotilde}(\lambda_x)
=
\frac{Z}{y_0^{\rtwotilde}}
<
\frac{Z}{y^{\rtwotilde}}
\leq1.
$$
Since $\Pi_{\rtwotilde}$ is continuous and strictly decreasing on
$[0,1/(\rtwotilde-1)]$, there is a unique
$\lambda_y\in[0,\lambda_x)$ with
$$
\Pi_{\rtwotilde}(\lambda_y)
=
\frac{Z}{y^{\rtwotilde}}.
$$
Note that the quantity $\lambda_x$  is exactly the one supplied by
Lemma~\ref{lem:higher-product-comparisons} for this value of $x$; in
particular $\lambda_x\leq1/(\rthreetilde-1)$ by Step~2 of its proof, and
we therefore also have $\lambda_y<\lambda_x\leq1/(\rthreetilde-1)$, which allows us to use Lemma~\ref{lem:higher-product-bound}. By  Lemma~\ref{lem:higher-product-bound},
\eqref{eq:higher-normalized-local}, and
\eqref{eq:higher-E},
\begin{equation}\label{eq:higher-large-necessary}
y^{\rthreetilde+(\rthreetilde-\rtwotilde)/(\rtwotilde-1)}\Pi_{\rthreetilde}(\lambda_y)
\leq
\frac{Y(\rthreetilde-\rtwotilde)sK(u-1)}{\rtwotilde}
E(u)^{(\rthreetilde-\rtwotilde)/(\rtwotilde-1)}.
\end{equation}

Substituting
$$
y^{\rtwotilde}
=
\frac{Z}{\Pi_{\rtwotilde}(\lambda_y)}
$$
into \eqref{eq:higher-large-necessary} gives
\begin{equation}\label{eq:higher-large-lambday}
Z^{\frac{\rthreetilde-1}{\rtwotilde-1}}
\frac{\Pi_{\rthreetilde}(\lambda_y)}
{\Pi_{\rtwotilde}(\lambda_y)^{\frac{\rthreetilde-1}{\rtwotilde-1}}}
\leq
W E(u)^{(\rthreetilde-\rtwotilde)/(\rtwotilde-1)} 
\end{equation}
because $y^{\rthreetilde+{\frac{\rthreetilde-1}{\rtwotilde-1}}-1}=(y^{\rtwotilde})^{\frac{\rthreetilde-1}{\rtwotilde-1}}$ and $
\frac1{\rtwotilde}
\left(
\rthreetilde+
\frac{\rthreetilde-\rtwotilde}{\rtwotilde-1}
\right)
=
\frac{\rthreetilde-1}{\rtwotilde-1},
$.

The function
$$
\lambda\longmapsto
\frac{\Pi_{\rthreetilde}(\lambda)}
{\Pi_{\rtwotilde}(\lambda)^{(\rthreetilde-1)/(\rtwotilde-1)}}
$$
is strictly decreasing.  This follows by differentiating its logarithm.
Indeed, the average of $h/(1-h\lambda)$ over $\rtwotilde\leq h\leq \rthreetilde-1$ is larger than
the corresponding average over $1\leq h\leq \rtwotilde-1$, which follows because $h\mapsto h/(1-h\lambda)$ is increasing.

Since
$\lambda_y<\lambda_x$, replacing $\lambda_y$ by $\lambda_x$ in
\eqref{eq:higher-large-lambday} strictly decreases its left-hand side,
so it suffices to contradict \eqref{eq:higher-large-lambday} by proving
\begin{equation}\label{eq:higher-large-target}
Z^{\frac{\rthreetilde-1}{\rtwotilde-1}}
\frac{\Pi_{\rthreetilde}(\lambda_x)}
{\Pi_{\rtwotilde}(\lambda_x)^{\frac{\rthreetilde-1}{\rtwotilde-1}}}
\geq
W\,E(u)^{{\frac{\rthreetilde-1}{\rtwotilde-1}}-1}.
\end{equation}

Using $y^{\rtwotilde}=Z/\Pi_{\rtwotilde}(\lambda_y)$ in
\eqref{eq:higher-large-necessary} and then replacing $\lambda_y$ by
$\lambda_x$ therefore makes its left side strictly smaller.

This is now a direct computation.  Since
$\Pi_{\rtwotilde}(\lambda_x)=Z/y_0^{\rtwotilde}$ and
$\rtwotilde {\frac{\rthreetilde-1}{\rtwotilde-1}}=\rthreetilde+{\frac{\rthreetilde-1}{\rtwotilde-1}}-1$, the left-hand side of
\eqref{eq:higher-large-target} equals
$y_0^{\rthreetilde+\frac{\rthreetilde-1}{\rtwotilde-1}-1}\Pi_{\rthreetilde}(\lambda_x)$.  Applying
\eqref{eq:higher-product-shift}, then
$\Pi_{\rthreetilde}(1/s)=YK/s^{\rthreetilde-1}$ and
$y_0=s/(1-x)$, gives
$$
\begin{aligned}
Z^{{\frac{\rthreetilde-1}{\rtwotilde-1}}}
\frac{\Pi_{\rthreetilde}(\lambda_x)}
{\Pi_{\rtwotilde}(\lambda_x)^{\frac{\rthreetilde-1}{\rtwotilde-1}}}
&=
y_0^{\rthreetilde+{\frac{\rthreetilde-1}{\rtwotilde-1}}-1}\Pi_{\rthreetilde}(\lambda_x)
\\
&\geq
\left(\frac{s}{1-x}\right)^{\rthreetilde+\frac{\rthreetilde-1}{\rtwotilde-1}-1}
\frac{YK}{s^{\rthreetilde-1}}
\bigl(1+(\rthreetilde-1)x\bigr)(1-x)^{\rthreetilde-1}
\bigl(1-(\rtwotilde-1)(\rthreetilde-1)x\bigr)^{{\frac{\rthreetilde-1}{\rtwotilde-1}}-1}
\\
&=
YKs^{{\frac{\rthreetilde-1}{\rtwotilde-1}}}
\frac{\bigl(1+(\rthreetilde-1)x\bigr)
\bigl(1-(\rtwotilde-1)(\rthreetilde-1)x\bigr)^{{\frac{\rthreetilde-1}{\rtwotilde-1}}-1}}
{(1-x)^{\frac{\rthreetilde-1}{\rtwotilde-1}}} .
\end{aligned}
$$
On the other hand, \eqref{eq:higher-large-x} gives
\begin{align*}
W\,E(u)^{{\frac{\rthreetilde-1}{\rtwotilde-1}}-1}
&=
\frac{YsK\bigl(1+(\rthreetilde-1)x\bigr)}{1-x}
\cdot
\frac{s^{{\frac{\rthreetilde-1}{\rtwotilde-1}}-1}\bigl(1-(\rtwotilde-1)(\rthreetilde-1)x\bigr)^{{\frac{\rthreetilde-1}{\rtwotilde-1}}-1}}
{(1-x)^{{\frac{\rthreetilde-1}{\rtwotilde-1}}-1}}\\
&=
YKs^{{\frac{\rthreetilde-1}{\rtwotilde-1}}}
\frac{\bigl(1+(\rthreetilde-1)x\bigr)
\bigl(1-(\rtwotilde-1)(\rthreetilde-1)x\bigr)^{{\frac{\rthreetilde-1}{\rtwotilde-1}}-1}}
{(1-x)^{{\frac{\rthreetilde-1}{\rtwotilde-1}}}} .
\end{align*}
The two expressions coincide, so
\eqref{eq:higher-large-target} holds.  Since replacing $\lambda_y$ by
$\lambda_x$ strictly decreased the left-hand side of
\eqref{eq:higher-large-lambday}, this contradicts
\eqref{eq:higher-large-lambday}.  Hence $y\geq y_0$, which is
\eqref{eq:higher-y-ge-y0} in the large-$u$ range.

Combining \eqref{eq:higher-central-identity} and
\eqref{eq:higher-y-ge-y0} gives
$$
M^{(\rtwotilde-1)/\rtwotilde}
\geq
1+\frac{(\rtwotilde-1)\widehat c}{\rtwotilde KD}.
$$
Since $M=(L+\widehat c)/(sK)$, this is exactly
\eqref{eq:higher-pointwise}.
\end{proof}

\begin{remark}
When $p=3$, we have $\rtwotilde=2$, and
\eqref{eq:higher-pointwise} becomes
$$
2KD-P
\geq
sK+\widehat c\left(\frac{s}{D}-1\right)
+\frac{s\widehat c^{\,2}}{4KD^2}.
$$
Thus the general pointwise bound recovers exactly the quadratic bound used
for the $(K_2,K_3,K_{q})$-profile.
\end{remark}

We now average the pointwise bound.  Its strict convexity converts the
local KKT inequalities into the global contradiction required by
Lemma~\ref{lem:higher-smooth-separation}.

\begin{proof}[Proof of Lemma~\ref{lem:higher-smooth-separation}]
Since $(s)_{\rtwotilde}=sK$, the prescribed means satisfy
\begin{equation}\label{eq:higher-separation-cancellation}
\frac{\rtwotilde K\overline D-\overline P}{\rtwotilde-1}
=\sigma sK\leq sK.
\end{equation}

Put
$$
\widehat c
=\frac{\rtwotilde c}
{(\rtwotilde-1)(\rthreetilde-\rtwotilde)}
$$
and
$$
\psi_c(D)
=sK\left(
1+\frac{(\rtwotilde-1)\widehat c}
{\rtwotilde KD}
\right)^{\rtwotilde/(\rtwotilde-1)}
-\widehat c.
$$
Direct differentiation shows that $\psi_c$ is strictly convex and decreasing on $D>0$.  Apply
Lemma~\ref{lem:higher-pointwise}.
Jensen's inequality then gives
$$
\frac{\rtwotilde K\overline D-\overline P}{\rtwotilde-1}
\geq\mathbb E\psi_c(D)
\geq\psi_c(\overline D)
\geq\psi_c(s).
$$
Here the last inequality uses
$\overline D=\sigma s(1-\alpha)\leq s$ and the fact that $\psi_c$ is
decreasing.
If
$$
x=\frac{(\rtwotilde-1)\widehat c}
{\rtwotilde sK}>0,
$$
strict Bernoulli gives
$$
\psi_c(s)
=sK(1+x)^{\rtwotilde/(\rtwotilde-1)}-\widehat c
>sK\left(1+\frac{\rtwotilde}{\rtwotilde-1}x\right)
-\widehat c
=sK.
$$
Thus the left side of \eqref{eq:higher-separation-cancellation} is
strictly greater than $sK$, contradicting that equation.  Notice that the
proof includes both endpoints $\alpha=0$ and $\alpha=1/s$; no limiting
argument is needed.
\end{proof}


\section*{Declaration on the use of AI}
Conversations with ChatGPT Pro 5.6 were instrumental in proving the technical lemmas. ChatGPT Pro 6 and Claude 5 Opus were helpful in writing the paper. All mathematical arguments were verified by the authors who take full responsibility for any mistakes. 

\section*{Acknowledgments}
We are grateful for conversations we had with Sasha Razborov, with colleagues at the Flag Algebras workshop at the American Institute of Mathematics in the Fall of 2025 as well as with colleagues at the Extremal Combinatorics semester at Simons Laufer Mathematical Sciences Institute in the Spring of 2025.

\end{document}